\documentclass[a4paper,11pt,twoside]{article}
\usepackage{pst-fill,pst-grad} 
\usepackage{textcomp}
\usepackage[titletoc]{appendix} 
\usepackage{titlesec}
\usepackage[english]{babel} 
\usepackage[utf8]{inputenc}  
\usepackage{graphicx} 
\usepackage{amsmath}
\usepackage{float} 
\usepackage{fancyhdr}
\usepackage[matrix,arrow,curve]{xy}
\usepackage{pstricks} 
\usepackage{amsmath,amsfonts,verbatim,afterpage,theorem,euscript,mathrsfs,amssymb}
\usepackage{amsfonts}
\usepackage{amssymb}
\usepackage{array}
\usepackage{dsfont}
\usepackage{hyperref}
\usepackage{authblk}

\newtheorem{Proposition}{Proposition}[section]

\newtheorem{Lemme}{Lemma}[section]
\newtheorem{Theoreme}{Theorem} 
\newtheorem{Corollaire}{Corollary}[section]
\newtheorem{Remarque}{Remark}[section]

\def \vf{\vec{f}} 
\def \vg{\vec{g}} 
\def \vu{\vec{u}}
\def \vn{\vec{\nabla}}
\def \vphi{\vec{\varphi}}
\def \P{\mathbb{P}}
\def \Rt{\mathbb{R}^{3}}

\def \finpv{\hfill $\blacksquare$  \newline }
\def \pv{{\bf{Proof.}}~}
\def \ds{\displaystyle}

\numberwithin{equation}{section}
\title{\bf  A turbulent study for a damped Navier-Stokes equation: turbulence and problems}
\author[1]{\small Diego Chamorro\footnote{\emph{Corresponding author}: diego.chamorro@univ-evry.fr}}
\author[2]{\small Oscar Jarr\'in\footnote{oscar.jarrin@udla.edu.ec}}
\affil[1]{\scriptsize LaMME, Univ Evry, CNRS, Université Paris-Saclay, 91025, Evry, France.}
\affil[2]{\scriptsize Escuela de Ciencias Físicas y Matemáticas, Universidad de Las Américas, Vía a Nayón, C.P.170124, Quito, Ecuador.}
\begin{document}
\maketitle
\begin{scriptsize}
\abstract{In this article we consider a damped version of the incompressible Navier-Stokes equations in the whole three-dimensional space with a divergence-free and time-independent external force. Within the framework of  a well-prepared force and with a particular choice of the damping parameter, when the Grashof numbers are large enough, we are able to prove some estimates from below and from above between the fluid characteristic velocity and the energy dissipation rate according to the Kolmogorov dissipation law. Precisely, our main contribution concerns the estimate from below which is not often studied in the existing literature. Moreover, we address some remarks which open the door to a deep discussion on the validity of this theory of turbulence.  
}\\[3mm]
\end{scriptsize}
\section{Introduction} 
The study of turbulence is an open and difficult problem in the field of fluid dynamics which has many implications in the mathematical, physical and technical world. In mathematics, different techniques have been used to try to solve this problem: indeed, from the stochastic point of view we can cite the works \cite{Chevillard}, \cite{Poursina} and \cite{Pozoroski} and from the deterministic point of view we can cite \cite{Chamorro}, \cite{DoerFoias}, \cite{FMRT1}, \cite{FMRT2} and \cite{Vinegron}.\\

In his celebrated 1941 theory, Andrey Kolmogorov \cite{Kolmogorov1} stated some simple laws in order to describe the turbulent regime and in the present article we will study the \emph{Kolmogorov dissipation law} in the deterministic framework. Let us first recall the global framework of the theory: we consider here the classical 3D Navier-Stokes equations 
\begin{equation}\label{Equation_NS_Intro}
\left\lbrace \begin{array}{lc}\vspace{3mm}
\partial_t\vu=\nu\Delta\vu-(\vu\cdot \vec{\nabla}) \vu-\vn p+\vf, \quad div(\vu)=0,\quad \nu>0,  &\\
\vu(0,\cdot)=\vu_0,&\\ 
\end{array}\right.
\end{equation}
where $\vu:[0,+\infty[\times \Rt\longrightarrow \Rt$ is the velocity of the fluid, $p:[0,+\infty[\times \Rt\longrightarrow\mathbb{R}$ is the pressure, $\nu >0$ is the fluid viscosity parameter and $\vf\in L^2(\Rt)\cap\dot{H}^{-1}(\Rt)$ is a given divergence-free, time-independent external force. Now, from a divergence-free initial data $\vu_0 \in L^2(\Rt)$, we can construct Leray solutions $\vu\in L^{\infty}_{loc}([0,+\infty[,L^2(\Rt)) \cap L^{2}_{loc}([0,+\infty[, \dot{H}^{1} (\Rt))$ which satisfy the energy inequality
\begin{equation}\label{Inegalite_Energie_Intro}
\|\vu(t, \cdot)\|_{L^2}^2+2\nu\int_{0}^t\|\vu(s, \cdot)\|_{\dot{H}^1}^2ds\leq \|\vu_0\|^2_{L^2}+2\int_{0}^{t} \langle \vf, \vu(s,\cdot)\rangle_{\dot{H}^{-1} \times \dot{H}^1} ds. 
\end{equation}
In this context, the time-independent external force $\vf$ will introduce a constant supply of energy (at some fixed length scale $\ell_0$ which represents the radius of a 3D ball where the force acts effectively over the fluid) and this permanent energy contribution is meant to lead the system to a turbulent state. Within this framework we define the averaged in time quantity
\begin{equation}\label{Def_U_Intro}
U= \left( \limsup_{T\to +\infty} \frac{1}{T} \int_{0}^{T} \| \vu(t,\cdot) \|^{2}_{L^2} \frac{dt}{\ell^{3}_{0}}\right)^{\frac{1}{2}},
\end{equation}
which represents the \emph{fluid characteristic velocity} and associated to this quantity we define the dimensionless \emph{Reynolds number} (see \cite{Constantin}, \cite{DoerFoias}) by 
\begin{equation}\label{Def_Re_Intro}
Re= \frac{U \ell_0}{\nu}. 
\end{equation} 
Finally, we introduce the \emph{energy dissipation rate} by the following averaged expression
\begin{equation}\label{Def_Epsilon_Intro}
\mathcal{E}= \nu \,    \limsup_{T\to +\infty} \frac{1}{T} \int_{0}^{T} \|\vu(t,\cdot)\|^{2}_{\dot{H}^1} \frac{dt}{\ell^{3}_{0}}.
\end{equation}
We recall that the turbulent regime of a fluid is usually characterized by the large values of the Reynolds numbers (\cite{Houel}, \cite{Tennekes}), then the Kolmogorov dissipation law \cite{Kolmogorov1} in the turbulent state reads as follows:
\begin{equation}\label{Kolmogorov_Law_Intro}
\mathcal{E}\sim \frac{U^3}{\ell_0} \qquad (Re\gg 1).
\end{equation}
This seemingly simple relationship raises some very deep issues. First note that the Reynolds number (\ref{Def_Re_Intro}) is built from the characteristic velocity $U$ which can be seen as an \emph{a posteriori} measurement and it will be interesting to replace the condition $Re\gg 1$ by another one which does not rely on $U$. Second, and far more annoying, we observe (as pointed out in \cite[Section 1.1.3]{Jarrin}) that since we consider the Navier-Stokes equations in the whole three-dimensional space, we can not ensure that the term $U$ given in (\ref{Def_U_Intro}) is a well-defined (bounded) quantity in the setting of Leray solutions and this last fact makes very hard to establish the relationship (\ref{Kolmogorov_Law_Intro}) rigorously. Note however that -by the energy inequality (\ref{Inegalite_Energie_Intro})- the quantity $\mathcal{E}$ given in (\ref{Def_Epsilon_Intro}) is well defined, see the Appendix \ref{Secc_AppendixA} for the details.\\

Let us mention that in a previous work \cite{Chamorro}, we introduced a suitable frequency truncation which allowed us to give a rigorous sense of these quantities. However, this modification appeared to be too strong: although the Reynolds number are large and that we have the relationship (\ref{Kolmogorov_Law_Intro}), we could proof that the fluid is not in a turbulent regime (see the details in \cite{Chamorro}).\\

To overcome some of these issues, we propose in this article the following modified Navier-Stokes equations:
\begin{equation}\label{Eq_damped_NS}
\left\lbrace \begin{array}{lc}\vspace{3mm}
\partial_t\vu=\nu\Delta\vu-(\vu\cdot \vec{\nabla}) \vu-\vn p+\vf-\beta \vu, \quad div(\vu)=0,\qquad (\beta,\, \nu>0),  &\\
\vu(0,\cdot)=\vu_0,&\\
\end{array}\right.
\end{equation}
where a damping term $-\beta \vu$ is added. This perturbation of the Navier-Stokes equation is a more relevant model than a brutal truncation in the Fourier variable: indeed, the parameter $0<\beta$ is known as the Rayleigh (or Ekman) friction coefficient and the term $-\beta \vu$ models the bottom friction in ocean models and is the main energy sink in large scale atmospheric models. See the book \cite[Chapter 5]{Pedlosky} for more details about this model. Of course when $\beta=0$ we recover the classical Navier-Stokes equations, but as we shall see the term $-\beta \vu$ will play a crucial role in our computations: in particular it will helps us to define correctly the quantity $U$ given in (\ref{Def_U_Intro}). \\

Next, we introduce the dimensionless \emph{Grashof number} (studied in \cite{DoerFoias} and \cite{Tennekes}) by the expression
\begin{equation}\label{Def_Gr_Intro}
Gr = \frac{\|\vf\|_{L^2} \ell^{\frac{3}{2}}_{0}}{\nu^2}.
\end{equation} 
Note that the Grashof number depends only on the external force $\vf$, its length scale $\ell_0$ and on the viscosity parameter $\nu$ but it is not related to the velocity field $\vu$. We thus shall see in the Corollary \ref{Coro_Relation_ReGr} below (under a particular choice of the external force $\vf$ and of the damping parameter $\beta$) that a large Grashof number ($Gr\gg 1$) implies that we have a large Reynolds number ($Re\gg 1$).\\

In the context of equation (\ref{Eq_damped_NS}) and with the help of the Grashof numbers defined in (\ref{Def_Gr_Intro}) we can prove (under some particular hypotheses over the external force $\vf$) the following version of the Kolmogorov dissipation law (\ref{Kolmogorov_Law_Intro}) 
\begin{equation}\label{Kolmogorov_Law_Intro1}
\mathcal{E}\sim \frac{U^3}{\ell_0} \qquad (Gr\gg 1),
\end{equation}
and the proof of the previous relationship constitutes the main result of this article (see Theorem \ref{Th_Kolmogorow_Law} below).\\

We continue by considering another modification of the classical Navier-Stokes equations (\ref{Equation_NS_Intro}) and for $\frac{3}{7} < \alpha < 3$ (the lower and upper bounds here are only technical, see the Remark \ref{Rem_LowerUpperBoundAlpha} below) we study the following fractional damped Navier-Stokes equations
\begin{equation}\label{Equation_FracNS_Intro}
\left\lbrace \begin{array}{lc}\vspace{3mm}
\partial_t\vu=-\nu(-\Delta)^{\frac{\alpha}{2}}\vu-(\vu\cdot \vec{\nabla}) \vu-\vn p+\vf-\beta \vu, \quad div(\vu)=0,\quad (\beta, \nu>0),  &\\
\vu(0,\cdot)=\vu_0,&\\
\end{array}\right.
\end{equation}
where the operator $(-\Delta)^{\frac{\alpha}{2}}$ can be defined in the Fourier level by its symbol $|\xi|^{\alpha}$. The fractional diffusion term $-\nu (-\Delta)^{\frac{\alpha}{2}} \vu $  has been successfully employed to model anomalous reaction-diffusion process in porous media models \cite{Meerschaert,Meerschaert2} and in computational turbulence models \cite[Chapter 13.2]{Pope}. In these last models, the term  $-\nu(-\Delta)^{\frac{\alpha}{2}}$ is used to  characterize anomalous viscous diffusion effects in turbulent fluids which are driven by the parameters $\alpha$ and $\nu$.\\

Our aim here is to establish in this setting a fractional version of the Kolmogorov dissipation law (\ref{Kolmogorov_Law_Intro1}) and to discuss the effects of the parameter $\frac{3}{7} < \alpha < 3$ in the turbulent regime. First, by a study based on dimensionality, we will obtain the following fractional dissipation law: 
\begin{equation}\label{Kolmogorov_FRAC_Law_Intro1}
\mathcal{E}\sim \frac{U^3}{\ell_0^{\alpha-1}} \qquad (Gr\gg 1),
\end{equation}
and we will see that in the case $\frac{3}{7} < \alpha<1$ the lower regularizing effect of the operator  $(-\Delta)^{\frac{\alpha}{2}}$ allows to consider a nicer external force while in the case $1\leq \alpha<3$ we need to consider a rougher force to establish the relationship (\ref{Kolmogorov_FRAC_Law_Intro1}). \\

Motivated by the previous study and noting that in the proof of this dissipation law the nonlinear term $(\vu\cdot\vn)\vu$ plays a minor role, we consider now the following damped Stokes equation
\begin{equation}\label{Eq_damped_Stokes}
\left\lbrace \begin{array}{lc}\vspace{3mm}
\partial_t\vu=\nu\Delta\vu-\vn p+\vf-\beta \vu, \quad div(\vu)=0,\qquad (\beta,\, \nu>0),  &\\
\vu(0,\cdot)=\vu_0.&\\
\end{array}\right.
\end{equation}
If $\beta=0$ we recover the usual Stokes equation which is a linearisation of the Navier-Stokes equation (\ref{Equation_NS_Intro}). This equation is usually meaningful in the case when the Reynolds number is small, however, if we add a time-independent external force $\vf$ we face to the similar problems when trying to define the quantity $U$ given in (\ref{Def_U_Intro}). Thus, by adding the term $-\beta \vu$ we easily ensure that this term $U$ is correctly defined. Now with the same external force used to study the damped Navier-Stokes system (\ref{Eq_damped_NS}) (see Section \ref{Secc_Main_theorem} below) and with $\beta\gg1$ we can establish the dissipation law (\ref{Kolmogorov_Law_Intro1}) for the damped Stokes equation (\ref{Eq_damped_Stokes}). \\

This fact raises some issues: on the one hand, we have a damped Stokes equation which (in principle) should not generate a turbulent behavior and, on the other hand, we have the relationship (\ref{Kolmogorov_Law_Intro1}) when the Grashof number (associated to the external force $\vf$) is large. This apparent contradiction may probably suggest that the Kolmogorov dissipation law (\ref{Kolmogorov_Law_Intro1}) is not an accurate model to describe turbulence in the whole space $\Rt$ and perhaps numerical simulations could helpful to a better understanding of the effect of the damping term coupled with the presence of the external force. \\

The outline of the article is the following. In Section \ref{Secc_Existence_Energy} we establish existence results as well as the energy inequalities for the damped Navier-Sokes equation (\ref{Eq_damped_NS}). In Section \ref{Secc_Main_theorem} we give the precise setting that allows us to deduce the Kolmogorov dissipation law (\ref{Kolmogorov_Law_Intro1}) and we state our main theorem (Theorem \ref{Th_Kolmogorow_Law}). Section \ref{Secc_Proof_MainTh} is devoted to the proof of this result. In Section \ref{Secc_remark_small_Grashof} we point out some important issues that are observed in our framework. In Section \ref{Secc_Stokes} we study the damped Stokes equation (\ref{Eq_damped_Stokes}) and finally, in Section \ref{Secc_FracNS}, we study the turbulence in the setting of the fractional equation (\ref{Kolmogorov_FRAC_Law_Intro1}).
\section{Existence of solutions, energy estimates and technical results}\label{Secc_Existence_Energy} 
We start with an existence result for the system (\ref{Eq_damped_NS}):
\begin{Theoreme}\label{Th1-Existence-NS-damped} Consider $\vu_0 \in L^{2}(\Rt)$ be a divergence-free initial data and let $\vf \in L^2(\Rt)\cap \dot{H}^{-1}(\Rt)$ be a divergence free and time independent  external force.  Then, for $0<\beta$, there exists a function $\vu\in L^{\infty}_{loc}([0,+\infty[,L^2(\Rt)) \cap L^{2}_{loc}([0,+\infty[, \dot{H}^{1} (\Rt))$ which is a weak solution of the equation (\ref{Eq_damped_NS}) that satisfies the energy inequality: 
\begin{equation}\label{Energy-inequality}
\begin{split}
\| \vu(t,\cdot)\|^{2}_{L^2}  \leq \,\, \| \vu_0 \|^{2}_{L^2} -  \, 2 \nu \, \int_{0}^{t}  \left\|\vu(s,\cdot)\right\|^{2}_{\dot{H}^{1}} ds + 2 \int_{0}^{t} \langle \vf, \vu(s,\cdot) \rangle_{L^2} \, ds -2 \beta \int_{0}^{t}\| \vu(s,\cdot)\|^{2}_{L^2} ds.
\end{split}
\end{equation} 
\end{Theoreme}
{\bf Proof.}  The proof of this theorem is rather classical and it essentially follows the Leray method for the classical Navier-Stokes equation. Let $\phi \in \mathcal{C}^{\infty}_{0}(\Rt)$ be a positive function such that $\ds{\int_{\Rt}\phi(x)dx=1}$. For $\delta>0$ we define $\phi_\delta(x)=\frac{1}{\delta^3}\phi \left( \frac{x}{\delta}\right)$.  Moreover, we denote by $h_{\nu t}(x)$ the usual heat kernel and with the help of the classical Leray projector $\mathbb{P}$ (defined for a suitable function $\vphi$ by $\mathbb{P}(\vphi)=\vphi+\vn \frac{1}{(-\Delta)}(\vn\cdot \vphi)$ and which is used here to get rid of the pressure $p$, see the book \cite{PGL} for the details), we will solve, due to the Duhamel formula, the following integral equation 
\begin{eqnarray}\label{NS_regularised_point_fix}\nonumber
\vu(t,x)&=&h_{\nu t}\ast \vu_0(x) +\int_{0}^{t}h_{\nu(t-s)}\ast \vf(x)ds-\int_{0}^{t}h_{\nu(t-s)}\ast(\P(([\phi_\delta\ast\vu]\cdot \vec{\nabla})\vu)(s,x)ds\\ 
& &-\beta \int_{0}^{t}h_{\nu(t-s)}\ast \vu(s,x)ds,
\end{eqnarray}
in the functional space $L^{\infty}([0,T],L^2(\Rt))\cap L^2([0,T],\dot{H}^1(\Rt))$ which is endowed, for some fixed $T>0$, with the norm $\Vert \cdot \Vert_T=\Vert \cdot \Vert_{L^{\infty}_{t}L^{2}_{x}}+\sqrt{\nu}\Vert \cdot \Vert_{L^{2}_{t}\dot{H}^{1}_{x}}$. We have then:
\begin{eqnarray*}
	\Vert \vu \Vert_T &\leq &\underbrace{\left\Vert h_{\nu t}\ast \vu_0 +\int_{0}^{t}h_{\nu(t-s)}\ast\vf(\cdot)ds\right\Vert_T}_{(A)}  + \underbrace{\left\Vert \int_{0}^{t}h_{\nu(t-s)}\ast(\P(([\phi_\delta\ast\vu]\cdot \vec{\nabla})\vu)(s,\cdot)ds\right\Vert_T}_{(B)} \\
	& & +\underbrace{\beta  \left\Vert \int_{0}^{t}h_{\nu(t-s)}\ast\vu(s,\cdot)ds \right\Vert_T}_{(C)}.
\end{eqnarray*}
The terms $(A)$ and $(B)$ are classical to estimate and by \cite{PGL}, Theorem $12.2$, page $352$,  we have $\Vert h_{\nu t}\ast \vu_0\Vert_{T}\leq c \Vert \vu_0 \Vert_{L^2}$ and  
\begin{equation}\label{estimate_force_th_existence}
\left\Vert \int_{0}^{t}h_{\nu(t-s)}\ast\vf(\cdot)ds \right\Vert_T \leq C(\frac{1}{\sqrt{\nu}}+\sqrt{T}) \Vert \vf \Vert_{L^{2}_{t}H^{-1}_{x}},
\end{equation} 
but since $\vf \in L^{2}(\Rt)\cap \dot{H}^{-1}(\Rt)$ we have: $\Vert \vf \Vert_{L^{2}_{t}H^{-1}_{x}}\leq \ds{\Vert \vf \Vert_{L^{2}_{t} \dot{H}^{-1}_{x} } \leq \sqrt{T} \Vert \vf \Vert_{L^{\infty}_{t} \dot{H}^{-1}_{x}} \leq \sqrt{T} \Vert \vf \Vert_{\dot{H}^{-1}}}$, and thus, for the term $(A)$ above we can write  
\begin{equation}\label{(1)}
\left\Vert h_{\nu t}\ast \vu_0 +\int_{0}^{t}h_{\nu(t-s)}\ast\vf(\cdot)ds\right\Vert_T  \leq  c\Vert \vu_0 \Vert_{L^2}+ C(\frac{1}{\sqrt{\nu}}+\sqrt{T}) \sqrt{T}\Vert \vf \Vert_{\dot{H}^{-1}}.
\end{equation} 
For the term $(B)$ we have (see \cite{PGL}, Theorem $12.2$ for the details):
\begin{equation}\label{(2)}
\left\Vert \int_{0}^{t}h_{\nu(t-s)}\ast(\P(([\phi_\delta\ast\vu]\cdot \vec{\nabla})\vu)(s,\cdot)ds\right\Vert_T \leq C (\frac{1}{\sqrt{\nu}}+1) \sqrt{T}\delta^{-\frac{3}{2}} \Vert \vu \Vert_T\,\Vert \vu \Vert_T.
\end{equation}
Finally, in order to study   the term $(C)$ we use again the estimate (\ref{estimate_force_th_existence}) (with $\vu$ instead of $\vf$) to write
\begin{equation}\label{(3)}
\begin{split}
\beta \left\Vert \int_{0}^{t} h_{\nu (t-s)} \ast \vu (s,\cdot) ds \right\Vert_{T}  \leq &\, \beta\,   C ( \frac{1}{\sqrt{\nu}} + \sqrt{T}) \Vert \vu \Vert_{L^{2}_{t} H^{-1}_{x}}  \leq C ( \frac{1}{\sqrt{\nu}} + \sqrt{T}) \Vert \vu \Vert_{L^{2}_{t} L^{2}_{x}} \\
\leq &\, \beta\,   C ( \frac{1}{\sqrt{\nu}} + \sqrt{T})  \sqrt{T} \Vert \vu \Vert_{L^{\infty}_{t} L^{2}_{x}} \leq \beta\,   C ( \frac{1}{\sqrt{\nu}} + \sqrt{T})  \sqrt{T} \Vert \vu \Vert_{T}. 
\end{split}
\end{equation}
Once we have inequalities (\ref{(1)}), (\ref{(2)}) and (\ref{(3)}) at our disposal,   for a  time $T>0$ small enough and for $\delta>0$, by the Banach contraction principle we obtain $\vu_\delta \in L^{\infty}([0,T],L^2(\Rt))\cap L^2([0,T],\dot{H}^1(\Rt))$ a local solution of  the equations (\ref{NS_regularised_point_fix}).\\

Now we will prove that this solution $\vu_\delta$ is global. We remark that the function $\vu_\delta$ satisfies the regularized equation 
\begin{equation}\label{NS-damped-regularized}
\partial_t\vu_\delta= \nu\Delta \vu_\delta-\P(([\phi_\delta\ast\vu_\delta]\cdot\vec{\nabla}) \vu_\delta)+\vf- \beta \vu_\delta,
\end{equation}
where all the terms belong to the space $L^2([0,T], H^{-1}(\Rt))$ and then we can write 
\begin{equation}\label{energ_eq_delta}
\begin{split}
	\frac{d}{dt} \Vert \vu_\delta(t,\cdot)\Vert^{2}_{L^2} = & \, 2\langle \partial_t \vu_\delta(t,\cdot),\vu_\delta(t,\cdot)\rangle_{H^{-1}\times H^1}\\
= &\, -2\nu \Vert \vu_\delta(t,\cdot)\Vert^{2}_{\dot{H}^1}+2 \langle \vf,\vu_\delta(t,\cdot)\rangle_{\dot{H}^{-1} \times \dot{H}^1}-2\beta \Vert \vu_{\delta} (t,\cdot) \Vert^{2}_{L^2}.
\end{split}
\end{equation}
But $-2\beta \Vert \vu_\delta(t,\cdot)\Vert^{2}_{L^2}$ is a negative quantity and we get
\begin{equation*}
\begin{split}
	\frac{d}{dt} \Vert \vu_\delta(t,\cdot)\Vert^{2}_{L^2}  \leq &\, -2\nu \Vert \vu_\delta(t,\cdot)\Vert^{2}_{\dot{H}^1} + 2\langle f,\vu_\delta(t,\cdot)\rangle_{\dot{H}^{-1}\times \dot{H}^1} \\
	\leq &\, -2\nu \Vert \vu_\delta (t,\cdot)\Vert^{2}_{\dot{H}^1} +\nu \Vert \vu_\delta(t,\cdot)\Vert^{2}_{\dot{H}^1}+\frac{1}{\nu}\Vert f \Vert^{2}_{\dot{H}^{-1}} \\
 \leq &\,  -\nu \Vert \vu_\delta (t,\cdot)\Vert^{2}_{\dot{H}^1}+\frac{1}{\nu}\Vert f \Vert^{2}_{\dot{H}^{-1}}.
\end{split}
\end{equation*}

Thereafter, we integrate on the interval of time $[0,t]$  and  we obtain the following control
\begin{equation}\label{control_unif_delta}
\Vert \vu_\delta (t,\cdot)\Vert^{2}_{L^2}+\nu\int_{0}^{t} \Vert \vu_\delta(s,\cdot)\Vert^{2}_{\dot{H}^1}ds \leq \left(\Vert\vu_0\Vert^{2}_{L^2}+\frac{t}{\nu}\Vert f \Vert^{2}_{\dot{H}^{-1}}\right), 
\end{equation} 
which allows us to extend the local solution $ \vu_\delta $ to the whole interval $[0,+\infty[$.\\
\\
We study now the convergence to a weak solution of the equations (\ref{Eq_damped_NS}). Indeed, by the  Rellich-Lions lemma (see \cite{PGL}, Theorem $12.1$) there exists a sequence of positive numbers $(\delta_n)_{n\in\mathbb{N}}$ and a function $\vu \in L^{2}_{loc}([0,+\infty[\times \Rt)$ such that the sequence $(\vu_{\delta_n})_{n\in \mathbb{N}}$  converges strongly to $\vu$ in $L^{2}_{loc}([0,+\infty[\times \Rt)$. Moreover,  this sequence converges  to $\vu$  in the weak$-*$ topology of the spaces  $L^{\infty}([0,T],L^2(\Rt))$ and $L^{2}([0,T],\dot{H}^1(\Rt))$ for all $T>0$. 
From these convergences we can deduce that the sequence $ \left(\P(([\phi_{\delta_n}\ast\vu_{\delta_n}]\cdot\vec{\nabla}) \vu_{\delta_n})\right)_{n\in \mathbb{N}}$ converges  to $\P((\vu\cdot\vec{\nabla}) \vu)$ in the weak$-*$ topology of the space  $(L^{2}_{t})_{loc}(H^{-\frac{3}{2}}_{x})$ and then, the limit $\vu$ is a weak  solution of equation  (\ref{Eq_damped_NS}). \\

Finally, in order to obtain  the energy inequality (\ref{Energy-inequality}), we get back to the identity (\ref{energ_eq_delta}) and we integrate (in the time variable) each term of this identity:
$$\Vert \vu_{\delta_n}(t,\cdot)\Vert^{2}_{L^2}+2\nu \int_{0}^{t} \Vert \vu_{\delta_n}(s,\cdot)\Vert^{2}_{\dot{H}^1}ds = \Vert \vu_0\Vert^{2}_{L^2} +2 \int_{0}^{t} \langle \vf, \vu_{\delta_n}(s,\cdot)\rangle_{L^2}ds -2\beta \int_{0}^{t}\Vert \vu_{\delta_n}(s,\cdot)\Vert^{2}_{L^2}ds,$$
from which we obtain the wished inequality (\ref{Energy-inequality}) by applying classical tools (see the book \cite{PGL}). Theorem \ref{Th1-Existence-NS-damped} is proven. \finpv 

With this theorem at our disposal, we will now show that the characteristic velocity $U$ given in the expression (\ref{Def_U_Intro}) is well defined. Indeed, we have:
\begin{Proposition}\label{Prop1}  Within the framework of Theorem \ref{Th1-Existence-NS-damped}, the solutions of the equation (\ref{Eq_damped_NS}) verify:
\begin{equation}\label{Estimation_Uf}
\limsup_{T \to +\infty} \frac{1}{T} \int_{0}^{T} \| \vu(t,\cdot) \|^{2}_{L^2} dt \leq \frac{1}{\nu\beta}\,\|\vf\|^{2}_{\dot{H}^{-1}}. 
\end{equation}
\end{Proposition}
{\bf Proof.} From the energy inequality (\ref{Energy-inequality})  and since $\vf \in \dot{H}^{-1}(\Rt)$,  we can  write
\begin{equation*}
\begin{split}
2\beta \int_{0}^{T} \Vert \vu(t,\cdot)\Vert^{2}_{L^2} dt + 2\nu \int_{0}^{T}\Vert \vu(t,\cdot)\Vert^{2}_{\dot{H}^{1}}dt \leq &\, \Vert \vu_0 \Vert^{2}_{L^2} + 2 \int_{0}^{T} \langle \vf , \vu(t,\cdot)\rangle_{\dot{H}^{-1} \times \dot{H}^{1}} dt  \\
\leq &\,  \Vert \vu_0 \Vert^{2}_{L^2} + 2 \int_{0}^{T} \Vert \vf \Vert_{\dot{H}^{-1}} \Vert \vu(t,\cdot)\Vert_{\dot{H}^{1}} dt \\
\leq &\,  \Vert \vu_0 \Vert^{2}_{L^2} + \nu  \int_{0}^{T}\Vert \vu(t,\cdot)\Vert^{2}_{\dot{H}^{1}}dt + \frac{T}{\nu} \Vert \vf \Vert^{2}_{\dot{H}^{-1}}.  
\end{split}
\end{equation*}
We thus obtain 
\[ \beta  \int_{0}^{T} \Vert \vu(t,\cdot)\Vert^{2}_{L^2} dt \leq  \Vert \vu_0 \Vert^{2}_{L^2}  + \frac{T}{\nu} \Vert \vf \Vert^{2}_{\dot{H}^{-1}}, \]
from which, we divide each term by $T$ and by letting $T\to +\infty$, we get the wished estimate (\ref{Estimation_Uf}).   \hfill$\blacksquare$\\

Note that due to the definition of $U$ given in (\ref{Def_U_Intro}), the inequality (\ref{Estimation_Uf}) is equivalent to the estimate $U\leq \frac{1}{\sqrt{\nu \beta}}\frac{\|\vf\|_{\dot{H}^{-1}}}{\ell_0^{\frac{3}{2}}}$, thus, as long as we have $\beta>0$ and $\vf \in \dot{H}^{-1}(\Rt)$ we can give a sense to the quantity $U$. From now on we will assume that we have 
\begin{equation}\label{ControlUFH_1}
1<U\leq \frac{1}{\sqrt{\nu \beta}}\frac{\|\vf\|_{\dot{H}^{-1}}}{\ell_0^{\frac{3}{2}}}.
\end{equation}

\section{The main result}\label{Secc_Main_theorem}
In order to present the main theorem of this article, we need to consider a particular external force $\vf$ that satisfies some conditions and we also need to be more precise with the damping parameter $\beta$. Let us start with the external force and for a fixed length scale $\ell_0>1$, we assume from now on that $\vf \in L^{2}(\Rt)\cap \dot{H}^{-1}(\Rt)$ verifies  the following frequency localization 
\begin{equation}\label{Loc_frec_force}
supp\left( \widehat{\vf} \right) \subset \left\{|\xi| < \frac{\mathfrak{c}}{\ell_0}\right\},
\end{equation}
where $0<\mathfrak{c}$ is a fixed and dimensionless constant. Only for technical reasons (see the Remark \ref{Remarque_ConditionF} below) we shall consider in this section $0<\nu<2$, and then,  we will also assume that we have 
\begin{equation}\label{Def_F1}
\|\vf\|_{L^2}\gg \ell_0^{\frac{5}{2}}\qquad \mbox{and}\qquad  \ell_0^{\frac{9}{8}}\|\vf\|_{\dot{H}^{-1}}\gg \sqrt{\nu}\|\vf\|_{L^2}^{\frac{7}{4}}.
\end{equation}
Note that all the previous conditions (\ref{Loc_frec_force}) and (\ref{Def_F1}) can be satisfied simultaneously by considering for example a function which is constant in the Fourier level over the ball $B(0,\frac{\mathfrak{c}}{\ell_0})$ and such that 
\begin{equation}\label{Condition_LinftyForce}
1\ll\ell_0^4\leq \|\widehat{f}\|_{L^\infty}\leq \ell_0^{\frac{13}{3}}.
\end{equation}

We define now the \emph{averaged external force} by the quantity 
\begin{equation}\label{Def_F}
F= \frac{\| \vf \|_{L^2}}{\ell^{\frac{3}{2}}_{0}},
\end{equation}
and we fix once and for all the parameter $\beta>0$ by the relationship
\begin{equation}\label{Def_Beta}
\beta=F^{\frac{3}{2}}.
\end{equation}
Now, with the fluid viscosity parameter $0<\nu<2$ we define the \emph{Grashof number} by the expression
\begin{equation}\label{Def_Gr}
Gr = \frac{F \ell^{3}_{0}}{\nu^2}= \frac{\|\vf\|_{L^2} \ell^{\frac{3}{2}}_{0}}{\nu^2}.
\end{equation} 
We recall that $Gr$ is a dimensionless quantity and we note that this definition is the same as the one given in the expression (\ref{Def_Gr_Intro}) above.\\

With these preliminaries, we can now state our main result. 
\begin{Theoreme}\label{Th_Kolmogorow_Law} 
Let $\vu_0 \in L^2(\Rt)$ be a divergence-free initial  data, let $\vf \in L^2(\Rt)\cap \dot{H}^{-1}(\Rt)$ be the divergence-free external force which verifies (\ref{Loc_frec_force}) and (\ref{Def_F1}) for some fixed length scale $1 \ll \ell_0$ and for $0<\nu<2$. Moreover, let $0<\beta$ be the damping parameter fixed in (\ref{Def_Beta}). Let $\vu \in L^{\infty}_{t}L^{2}_{x}\cap (L^{2}_{t})_{loc}\dot{H}^{1}_{x}$ be a solution of the equation (\ref{Eq_damped_NS}) associated to $\vu_0, \vf$ and $\beta$.\\

\noindent If the Grashof number is large enough (\emph{i.e.} $Gr\gg 1$), then we have the estimates:
\begin{equation}\label{Estimations_TheoremePrincipal}
\frac{1}{2}\frac{U^3}{\ell_0}\leq \mathcal{E}\leq\left(\frac{101}{10}\right)\frac{U^3}{\ell_0},
\end{equation}
where the quantities $U$ and $\mathcal{E}$ are defined in (\ref{Def_U_Intro}) and (\ref{Def_Epsilon_Intro}) respectively. 
\end{Theoreme}
Some remarks are in order here. First note that our approach is based on the presence of a particular external force $\vf$ whose properties rely heavily on the length scale $\ell_0$ which is assumed to be large. Since we have $\|\vf\|_{L^2}\gg \ell_0^{\frac{5}{2}}$ by the condition (\ref{Def_F1}), the Grashof number $Gr$ given in (\ref{Def_Gr}) will be also large, thus, as the region where the force is effective is consequent (which is precisely the notion of $\ell_0$) and as the energy input is big, it is quite reasonable to obtain the dissipation law (\ref{Estimations_TheoremePrincipal}) as predicted by the general theory.  We also observe that if we have $\|\vf\|_{L^2}\gg \ell_0^{\frac{5}{2}}$, then the quantity $F$ given in (\ref{Def_F}) is also large and so is the parameter $\beta=F^{\frac{3}{2}}$ defined in (\ref{Def_Beta}): we thus obtain the wished dissipation law with an important damping term. 

\medskip

On the other hand, it is interesting to remark that with a well-prepared external force and with some very particular values of the  parameters in our model we can \emph{also} obtain similar estimates as in (\ref{Estimations_TheoremePrincipal}) even in the \emph{non turbulent case} when  $Gr \sim 1$ and $\ell_0=1$.  In Section \ref{Secc_remark_small_Grashof} below we address all the details. This   remark yields us to  question, on the one hand,  the relevance of the damping parameter $\beta$ in our study and, on the other hand, the validity of this turbulence theory. 

\medskip  

For our calculation to work we need a large damping parameter $\beta$ and perhaps some interference between the damping term and the expected dissipation law can be made explicit, but this point is not completely clear for us and perhaps numerical simulations can bring an interesting insight to the behavior of the fluid in this case. Let us make clear here that we do not claim any optimality in our approach and perhaps a different perturbation of the Navier-Stokes equation can show more explicitly how to obtain the Kolmogorov dissipation law (\ref{Kolmogorov_Law_Intro1}).

\medskip

Finally, in the Appendix \ref{Secc_AppendixC} we study the Kolmogorov dissipation  law (\ref{Kolmogorov_Law_Intro1}) with another different choice for the damping parameter $\beta$ (see the expression (\ref{Def_Beta_2} for a precise definition) which does not depend on the quantity $\Vert \vf \Vert_{L^2}$, and consequently, in this case the damping parameter is not necessary large as in (\ref{Def_Beta}). However, it is worth emphasizing that this different choice for $\beta$ does not allow us derive  (to the best of our knowledge) an estimate from bellow as in the estimates (\ref{Estimations_TheoremePrincipal}).  See Theorem \ref{Th_Kolmogorow_Law_Damping2} in the Appendix \ref{Secc_AppendixC} for all the details. 

\section{Proof of Theorem \ref{Th_Kolmogorow_Law}}\label{Secc_Proof_MainTh}
Note that in this framework all the quantities $\nu$, $\ell_0$, $F$ and $U$ are \emph{finite} and we need to establish some relationships between them in order to perform our computations. In this sense we have the following result.
\begin{Lemme}\label{Lem_EstimatesFUF}
Under the hypothesis of the Theorem \ref{Th_Kolmogorow_Law} we have the controls
\begin{equation}\label{Estimation_FUF0}
1\ll F\leq U\leq F^{\frac{3}{2}} \ell_0.
\end{equation}
\end{Lemme}
{\bf Proof.}  By the energy inequality (\ref{Energy-inequality})  we can write
$$2\beta \int_{0}^{T} \| \vu(t,\cdot) \|^{2}_{L^2} dt + \| \vu(T,\cdot) \|^{2}_{L^2} \leq \| \vu_0 \|^{2}_{L^2} + 2 \int_{0}^{T} \int_{\Rt} \vf\cdot \vu \, dx \, dt - 2 \nu \int_{0}^{T}\| \vu(t,\cdot)\|^{2}_{\dot{H}^{1}} dt,$$
from which we deduce the estimate
\begin{equation*}
\begin{split}
2\beta \int_{0}^{T} \| \vu(t,\cdot) \|^{2}_{L^2} dt \leq &\,  \| \vu_0 \|^{2}_{L^2} + 2 \int_{0}^{T} \int_{\Rt} \vf\cdot \vu \, dx \, dt \\
\leq &\, \| \vu_0 \|^{2}_{L^2} + 2 \int_{0}^{T}\| \vf \|_{L^2}\, \| \vu(t,\cdot) \|_{L^2} dt. \end{split}
\end{equation*}
Applying the Young inequalities with some parameter $\gamma>0$ in the integral above (recall that $\beta$ is given by (\ref{Def_Beta})), we obtain
$$2\beta \int_{0}^{T} \| \vu(t,\cdot) \|^{2}_{L^2} dt \leq  \| \vu_0 \|^{2}_{L^2} + 2 \int_{0}^{T}\frac{1}{\gamma \beta}\|\vf\|_{L^2}^2 +\gamma\beta\|\vu(t,\cdot)\|_{L^2}^2 dt. $$
which can be rewritten in the following manner:
$$(2-\gamma)\beta\int_{0}^{T} \|\vu(t,\cdot)\|^{2}_{L^2} dt\leq \|\vu_0\|^{2}_{L^2} + T\frac{1}{\gamma \beta}\|\vf\|_{L^2}^2$$
Dividing each term above by $\ell^{3}_{0}$ and by $T$ and taking the limit when $T\to +\infty$, we thus have
$$(2-\gamma)\beta\limsup_{T\to +\infty}\frac{1}{T} \int_{0}^{T} \|\vu(t,\cdot)\|^{2}_{L^2} \frac{dt}{\ell_0^3}\leq \frac{1}{\gamma \beta}\frac{\|\vf\|_{L^2}^2}{\ell_0^3},$$
now, by the definition of the quantities $U$ and $F$ given in (\ref{Def_U_Intro}) and (\ref{Def_F}) we have
\begin{equation}\label{Equation_GeneriqueUFU}
\gamma(2-\gamma)\beta^2 \, U^2 \leq F^2.
\end{equation}
We define now $\gamma$ such that $\gamma(2-\gamma)=\frac{1}{U^3}<1$ (which is the case since we assumed in (\ref{ControlUFH_1}) that $U>1$), thus with the fact that $\beta^2=F^3$ (see (\ref{Def_Beta})), the estimates above becomes
$$F\leq U.$$
If we consider now $\gamma$ such that $\gamma(2-\gamma)=\frac{1}{\ell_0^2 F^4}<1$
(Note that $\frac{1}{\ell_0^2 F^4}<1$ is equivalent to $\ell_0<\|\vf\|_{L^2}$ which is a consequence from the first  condition in (\ref{Def_F1}) and the fact that $\ell_0>1$), from (\ref{Equation_GeneriqueUFU}) we easily obtain $U\leq \ell_0 F^{\frac{3}{2}}$ and the Lemma \ref{Lem_EstimatesFUF} is proven. \hfill $\blacksquare$
\begin{Remarque}\label{Remarque_ConditionF}
Note that the condition $F\leq U$ that we have just proved is compatible with the constraint $U\leq \frac{1}{\sqrt{\nu \beta}}\frac{\|\vf\|_{\dot{H}^{-1}}}{\ell_0^{\frac{3}{2}}}$ given in (\ref{ControlUFH_1}) as long as we have the second condition over the force stated in (\ref{Def_F1}) with $0<\nu<2$. In this case we do have the estimates $\frac{\|\vf\|_{L^2}}{\ell_0^{\frac{3}{2}}}=F\leq U\leq \frac{1}{\sqrt{2 \beta}}\frac{\|\vf\|_{\dot{H}^{-1}}}{\ell_0^{\frac{3}{2}}}  \leq \frac{1}{\sqrt{\nu  \beta}}\frac{\|\vf\|_{\dot{H}^{-1}}}{\ell_0^{\frac{3}{2}}}$.
\end{Remarque}
From the point of view of the Grashof and Reynold numbers we have the following consequence of this result:
\begin{Corollaire}\label{Coro_Relation_ReGr}
From the lower estimate in (\ref{Estimation_FUF0}) and with the definition of $Re$ and $Gr$ given in (\ref{Def_Re_Intro}) and (\ref{Def_Gr}) respectively, we obtain 
$$Gr\leq \frac{\ell_0^2}{\nu}Re,$$
and thus if $Gr\gg1$ we also have that $Re\gg 1$.\\
\end{Corollaire}
With all these estimates at hand, we can now study the proof of the Theorem \ref{Th_Kolmogorow_Law}.
\begin{itemize}
\item {\bf Lower estimate.} By the energy inequality (\ref{Energy-inequality}) we have this time
$$2\beta \int_{0}^{T} \| \vu(t,\cdot) \|^{2}_{L^2} dt + \| \vu(T,\cdot) \|^{2}_{L^2} \leq \| \vu_0 \|^{2}_{L^2} + 2 \int_{0}^{T} \int_{\Rt} \vf\cdot \vu \, dx \, dt + 2 \nu \int_{0}^{T}\| \vu(t,\cdot)\|^{2}_{\dot{H}^{1}} dt,$$
from which we deduce
$$2\beta \int_{0}^{T} \| \vu(t,\cdot) \|^{2}_{L^2} dt  \leq \| \vu_0 \|^{2}_{L^2} + 2 \int_{0}^{T} \|\vf\|_{L^2}\|\vu(t,\cdot)\|_{L^2} \, dt + 2 \nu \int_{0}^{T}\| \vu(t,\cdot)\|^{2}_{\dot{H}^{1}} dt,$$
and by the Young inequalities we have
$$2\beta \int_{0}^{T} \| \vu(t,\cdot) \|^{2}_{L^2} dt  \leq \| \vu_0 \|^{2}_{L^2} + \int_{0}^{T} \|\vf\|_{L^2}^2+\|\vu(t,\cdot)\|_{L^2}^2 \, dt + 2 \nu \int_{0}^{T}\|\vu(t,\cdot)\|^{2}_{\dot{H}^{1}} dt,$$
dividing the previous estimate by $\ell^{3}_{0}$ and $T$, taking the limit when $T\to +\infty$ we obtain
\begin{eqnarray*}
2\beta \limsup_{T\to +\infty}\frac{1}{T} \int_{0}^{T} \| \vu(t,\cdot) \|^{2}_{L^2} \frac{dt}{\ell_0^3} &\leq& \frac{\|\vf\|_{L^2}^2}{\ell_0^3}+ \limsup_{T\to +\infty}\frac{1}{T} \int_{0}^{T} \|\vu(t,\cdot)\|_{L^2}^2 \, \frac{dt}{\ell_0^3}\\
&& + 2 \nu \limsup_{T\to +\infty}\frac{1}{T}\int_{0}^{T}\|\vu(t,\cdot)\|^{2}_{\dot{H}^{1}} \frac{dt}{\ell_0^3}.\end{eqnarray*}
Recalling at this point the definition of the quantities $U$ given in (\ref{Def_U_Intro}), $F$ given in (\ref{Def_F}) and $\mathcal{E}$ given in (\ref{Def_Epsilon_Intro}), we obtain
$$2\beta U^2\leq F^2+ U^2+2\mathcal{E},$$
from which we deduce the estimate
$$\left(\beta-\frac{1}{2}\right)U^2-\frac{1}{2}F^2\leq \mathcal{E}.$$
Now, since by the Lemma \ref{Lem_EstimatesFUF} we have the estimate $F^2\leq U^2$, we obtain
\begin{equation}\label{Estim01}
\left(\beta-1\right)U^2\leq\left(\beta-\frac{1}{2}\right)U^2-\frac{1}{2}F^2\leq \mathcal{E}.
\end{equation}
We thus have, 
$$\left(\frac{\beta \ell_0}{U}-\frac{ \ell_0}{U}\right)\frac{U^3}{\ell_0}\leq \mathcal{E}.$$
Now, by the definition  (\ref{Def_Beta}) we have $\beta=F^{\frac{3}{2}}$ and since by the upper bound of the Lemma \ref{Lem_EstimatesFUF} we have $1\leq \frac{F^{\frac{3}{2}}\ell_0}{U}$, we can thus write
$$\left(1-\frac{ \ell_0}{U}\right)\frac{U^3}{\ell_0}\leq \left(\frac{\beta \ell_0}{U}-\frac{ \ell_0}{U}\right)\frac{U^3}{\ell_0}\leq \mathcal{E}.$$
Noting now that by the lower bound of the Lemma \ref{Lem_EstimatesFUF} we have $-\frac{1}{F}\leq -\frac{1}{U}$, we have
$$\left(1-\frac{\ell_0}{F}\right)\frac{U^3}{\ell_0}\leq \mathcal{E}.$$
We use now the definition of the Grashof numbers given in (\ref{Def_Gr}) to obtain the estimate
\begin{equation}\label{Estimate_Kolmogorov_below}
\left(1-\frac{\ell_0^4}{\nu^2 Gr}\right)\frac{U^3}{\ell_0}\leq \mathcal{E}.
\end{equation}
At this point we remark that we have $\frac{\ell_0^4}{\nu^2 Gr}\leq \frac{1}{2}$. Indeed, this control is equivalent by (\ref{Def_Gr}) to $\frac{\ell_0^4}{\nu^2}\leq \frac{Gr}{2}=\frac{1}{2}\frac{\|\vf\|_{L^2}\ell_0^{\frac{3}{2}}}{\nu^2}$ which can be rewritten as $\ell_0^\frac{5}{2}\leq\frac{1}{2}\|\vf\|_{L^2}$. We only need to remark that by (\ref{Def_F1}) we have $\ell_0^\frac{5}{2}\ll\|\vf\|_{L^2}$ to obtain the claimed control. Whit this remark we have:
$$\frac{1}{2}\frac{U^3}{\ell_0}\leq \mathcal{E},$$
which is the lower estimate in (\ref{Estimations_TheoremePrincipal}).\\
\item {\bf Upper estimate.} 
By the energy inequality (\ref{Energy-inequality}) and by the Cauchy-Schwarz inequality (first in the spatial variable and then in the temporal variable) we can write 
\begin{equation*}
\begin{split}
2 \nu \int_{0}^{T} \| \vn\otimes \vu(t,\cdot) \|^{2}_{L^2}\, dt \leq &\,  \| \vu_0 \|^{2}_{L^2} + 2 \int_{0}^{T} \int_{\Rt} \vf(x)\cdot \vu(t,x) \, dx\, dt  \\
\leq &\, \| \vu_0 \|^{2}_{L^2} + 2 \int_{0}^{T} \| \vf \|_{L^2}\, \| \vu(t,\cdot)\|_{L^2}\, dt.
\end{split}
\end{equation*}
Now, by the Young inequalities (with $\rho>0$), we obtain
$$2 \nu \int_{0}^{T} \| \vn\otimes \vu(t,\cdot) \|^{2}_{L^2}\, dt \leq\| \vu_0 \|^{2}_{L^2} + \int_{0}^{T} \frac{1}{\rho}\| \vf \|_{L^2}^2 +\rho\| \vu(t,\cdot)\|_{L^2}^2\, dt.$$

Dividing all this estimate by  $\ell^{3}_{0}$ and by $T$ and taking the limit $T \to +\infty$ we thus have
$$2 \nu \limsup_{T\to+\infty}\frac{1}{T}\int_{0}^{T} \| \vn\otimes \vu(t,\cdot) \|^{2}_{L^2}\, \frac{dt}{\ell_0^3} \leq  \frac{1}{\rho}\frac{\|\vf \|_{L^2}^2}{\ell_0^3} + \rho\limsup_{T\to+\infty}\frac{1}{T}\int_{0}^{T} \| \vu(t,\cdot)\|_{L^2}^2\, \frac{dt}{\ell_0^3}.$$
Now, by definition of the quantities $\mathcal{E}$, $F$ and $U$ given in (\ref{Def_Epsilon_Intro}), (\ref{Def_F}) and (\ref{Def_U_Intro}) respectively, we obtain
\begin{equation}\label{Estim02}
2\mathcal{E}\leq  \frac{1}{\rho}F^2+\rho U^2.
\end{equation}
If we set $\rho=\frac{\ell_0}{10 U}$ we have
$$\mathcal{E}\leq  \frac{10 U}{\ell_0}F^2+\frac{\ell_0}{10}U.$$
By the lower bound of the Lemma \ref{Lem_EstimatesFUF} we have $F^2\leq U^2$ and we can write
$$\mathcal{E}\leq  \frac{10 U^3}{\ell_0}+\frac{\ell_0}{10}U=\frac{U^3}{\ell_0}\left(10+\frac{\ell_0^2}{10U^2}\right)\leq \frac{U^3}{\ell_0}\left(10+\frac{\ell_0^2}{10F^2}\right).$$
Using the definition of the Grashof numbers given in (\ref{Def_Gr}) we finally obtain
\begin{equation}\label{Estimate_Kolmogorov_above}
\mathcal{E}\leq\frac{U^3}{\ell_0}\left(10+\frac{\ell_0^8}{10\nu^4 Gr^2}\right).
\end{equation}
$$\mathcal{E}\leq\frac{U^3}{\ell_0}\left(10+\frac{\ell_0^8}{10\nu^4 Gr^2}\right).$$
We observe now that we have $\frac{\ell_0^8}{\nu^4 Gr^2}\leq 1$, indeed, this is equivalent to $\frac{\ell_0^4}{\nu^2}\leq  Gr=\frac{\|\vf\|_{L^2}\ell_0^\frac{3}{2}}{\nu^2}$ and we thus obtain the condition $\ell_0^{\frac{5}{2}}\leq \|\vf\|_{L^2}$, which is the case since we do have $\|\vf\|_{L^2}\gg \ell_0^{\frac{5}{2}}$ by (\ref{Def_F1}).\\
We can finally write 
$$\mathcal{E}\leq\frac{U^3}{\ell_0}\left(\frac{101}{10}\right),$$
which is the upper estimate of (\ref{Estimations_TheoremePrincipal}).\\
\end{itemize}
With these two inequalities, the proof of Theorem \ref{Th_Kolmogorow_Law} is finished.\hfill $\blacksquare$
\section{Kolmogorov type estimates for small Grashof numbers}\label{Secc_remark_small_Grashof}
In this section, we show that we can obtain similar estimates as in  (\ref{Estimations_TheoremePrincipal}) when $Gr \sim 1$. For this  we will set the parameters in our model as follows. First, we set the length scale   and the viscosity parameter  as 
\begin{equation}\label{ell_0_nu}
\ell_0=1 \quad \text{and}\quad  \nu = \sqrt{2}.
\end{equation} 
Then, we shall consider a particular external force $\vf$ defined in the Fourier level as 
\begin{equation}\label{Def_particular_force}
\widehat{\vf}(\xi)=  2^{\frac{5}{2}}\, \mathds{1}_{\vert \xi \vert < \frac{1}{2}} (\xi),
\end{equation}
where  by the Plancherel identity we have 
\begin{equation}\label{L2_norm}
\Vert \vf \Vert_{L^2}  = 2 \left( \frac{4 \pi}{3}\right)^{\frac{1}{2}}.
\end{equation}
We thus have $F= \frac{\Vert \vf \Vert_{L^2}}{\ell^{\frac{3}{2}}_{0}}=\Vert \vf \Vert_{L^2} \simeq 4$ and as in the last section  we set  $\beta = F^{\frac{3}{2}}$. Moreover,  by (\ref{Def_Gr}) we have   $$Gr = \frac{\Vert \vf \Vert_{L^2} \ell^{\frac{3}{2}}_{0}}{\nu^2} = \left( \frac{4 \pi}{3}\right)^{\frac{1}{2}} \simeq 2.$$ 
This particular external force also verifies the  following  estimates, which are similar to the ones given in (\ref{Def_F1})
\begin{equation}\label{Def_F1_particular}
\Vert \vf \Vert_{L^2} \geq  1  \quad \text{and} \quad  \Vert \vf \Vert_{\dot{H}^{-1}} \simeq \sqrt{\nu} \Vert \Vert \vf \Vert^{\frac{7}{4}}_{L^2}.  
\end{equation}

 Indeed, the first estimate is a direct consequence of the first identity in (\ref{ell_0_nu}) and the identity (\ref{L2_norm}).  For the second condition, on the one hand,  by  a simple computation we have $\Vert \vf \Vert_{\dot{H}^{-1}}=4 (4\pi)^{\frac{1}{2}} \simeq 14.16$. On the other hand, by the second identity in (\ref{ell_0_nu}) and by the identity (\ref{L2_norm}) we have $\sqrt{\nu} \Vert \vf \Vert^{\frac{7}{2}}_{L^2} = 4 \left( \frac{4 \pi}{3}\right)^{\frac{7}{8}}\simeq 14.16$.\\
 
\medskip
 
 Now, by following the same ideas in the proof of Theorem \ref{Th_Kolmogorow_Law} we will quickly  show that we are able to  obtain  the estimates $\mathcal{E} \sim \frac{U^3}{\ell_0}$ even if $Gr \simeq 2$ and $\ell_0=1$.
 
 \medskip
 
First we remark that the estimate (\ref{Estimation_FUF0}) proven in Lemma \ref{Lem_EstimatesFUF} above is still valid in this particular setting and since $\ell_0=1$ it writes down as 
\begin{equation}\label{Estimation_FUF0_particular} 
1 < F \leq U \leq F^{\frac{3}{2}}.
\end{equation}
Moreover, as pointed out in Remark \ref{Remarque_ConditionF}, we shall remark that  the condition $F \leq U$ is still compatible with the constraint $U\leq \frac{1}{\sqrt{\nu \beta}} \|\vf\|_{\dot{H}^{-1}} $ (with $\ell_0=1$) given in (\ref{ControlUFH_1}). Precisely, in this case we have $F \simeq U \simeq \frac{1}{\sqrt{\nu \beta}}  \Vert \vf \Vert_{\dot{H}^{-1}}$ which is equivalent to  the second condition in (\ref{Def_F1_particular}).  

\medskip

With the estimates (\ref{Estimation_FUF0_particular}) at our disposal,  we get the estimate from below given in (\ref{Estimate_Kolmogorov_below}) and   since $\ell_0=1$ we shall write this estimate  as 
$$ \left(1-\frac{1}{\nu^2 Gr}\right)\frac{U^3}{\ell_0}\leq \mathcal{E}.$$
But, since  $Gr \simeq 2$ and $\nu=\sqrt{2}$ we have $\frac{1}{\nu^2 Gr} \simeq \frac{1}{4}$ and we thus obtain the estimate from below
$$ \frac{3}{4}  \frac{U^3}{\ell_0}\leq \mathcal{E}.$$
Always by the estimates (\ref{Estimation_FUF0_particular}), we also get the estimate from above given in (\ref{Estimate_Kolmogorov_above}), which in this particular case  ($\ell_0=1$, $\nu=\sqrt{2}$ and $Gr \simeq 2$)  writes down  as:
$$  \mathcal{E}\leq\frac{U^3}{\ell_0}\left(10+\frac{1}{10\nu^4 Gr^2}\right) \simeq \frac{1601}{160} U^3.$$

\medskip 

Summarizing, by setting the damping parameter $\beta$ as in (\ref{Def_Beta}), on the one hand with the particular external force $\vf$ defined in (\ref{Def_particular_force}) and with the particular values of the parameters in our model given in (\ref{ell_0_nu}), we  obtain the Kolmogorov type estimates   $\mathcal{E} \sim \frac{U^3}{\ell_0}$ in the non turbulent regime $Gr \sim 1$ and $\ell_0=1$.  On the other hand,  when we  consider  the particular force given in the example (\ref{Condition_LinftyForce}), and with $\ell_0\gg 1$, we also obtain the estimates $\mathcal{E} \sim \frac{U^3}{\ell_0}$ in the turbulent regime $Gr \gg 1$.  

\medskip 

We thus observe that by carefully setting the external force and the  parameters of our model we  are able to obtain  Kolmogorov type estimates in both the turbulent and the non turbulent cases. As mentioned in Section \ref{Secc_Main_theorem}, it is thus quite natural to deeply question the significance of the parameter $\beta$ in our study  as well as the validity of this theory of the turbulence.

\section{Turbulence in the Stokes equations}\label{Secc_Stokes}
As pointed out in the introduction, we may observe in the proof of the Theorem \ref{Th_Kolmogorow_Law} that the study of the Kolmogorov dissipation law in the deterministic framework of the damped Navier-Stokes equations
is strongly based on the energy inequality. On the other hand, in this energy inequality we see that the nonlinear effects of the transport term $(\vu \cdot \vn ) \cdot \vu$ are neglected since  this term vanishes in the computations: as $div(\vu)=0$ we formally have the identity $\displaystyle{\int_{\Rt} (\vu \cdot \vec{\nabla}) \vu \cdot \vu \, dx=0}$.   This remark yields a deeper discussion on how the nonlinear transport term intervenes in the deterministic study of the Kolmogorov dissipation law, and further, on the robustness of  this theory of  turbulence. 

\medskip

To start our discussion,  we shall consider here the damped Stokes  equations, which essentially writes down as the damped Navier-Stokes equation  without the nonlinear transport term:
\begin{equation}\label{Eq_damped_Stokes1}
\left\lbrace \begin{array}{lc}\vspace{3mm}
\partial_t\vu=\nu\Delta\vu-\vn p+\vf-\beta \vu, \quad div(\vu)=0,\qquad (\beta,\, \nu>0),  &\\
\vu(0,\cdot)=\vu_0.&\\
\end{array}\right. 
\end{equation}
As we have a linear model, the global the well-posedness of this equation is a well-known fact and we have:
\begin{Proposition} Let $\vu_0 \in L^2(\Rt)$ be a divergence-free initial data and let $\vf \in L^2(\Rt)\cap\dot{H}^{-1}(\Rt)$ be a divergence-free external force. Then, for $0<\beta$ there exists a function 
$\vu \in L^{\infty}_{t}L^{2}_{x}\cap (L^{2}_{t})_{loc}\dot{H}^{1}_{x}$ which is the unique strong solution of the equation (\ref{Eq_damped_Stokes1}). Moreover, the following energy equality holds:
\begin{equation}\label{Energy-equality}
\begin{split}
\| \vu(t,\cdot)\|^{2}_{L^2}  = \,\, \| \vu_0 \|^{2}_{L^2} -  \, 2 \nu \, \int_{0}^{t}  \left\|\vu(s,\cdot)\right\|^{2}_{\dot{H}^{1}} ds + 2 \int_{0}^{t} \int_{\Rt} \vf\cdot \vu(s,\cdot) dx \, ds -2 \beta \int_{0}^{t}\| \vu(s,\cdot)\|^{2}_{L^2} ds.
\end{split}
\end{equation}   
\end{Proposition}
The energy equality (\ref{Energy-equality}) also yields  the Lemma \ref{Lem_EstimatesFUF} which is the key tool in the proof of Theorem \ref{Th_Kolmogorow_Law}. Consequently, we can follow the same estimates performed in this proof to obtain   analogous  lower and upper bounds on the energy dissipation rate according to the Kolmogorov dissipation law:
\begin{Theoreme} Under the same hypothesis of the Theorem \ref{Th_Kolmogorow_Law}, let  $\vu \in L^{\infty}_{t}L^{2}_{x}\cap (L^{2}_{t})_{loc}\dot{H}^{1}_{x}$ be the solution of the equation (\ref{Eq_damped_Stokes}) from which we define the quantities $U$ and $\mathcal{E}$ given in (\ref{Def_U_Intro}) and (\ref{Def_Epsilon_Intro}) respectively. Then, we have the estimates 
$$\frac{1}{2}\frac{U^3}{\ell_0}\leq \mathcal{E}\leq\left(\frac{101}{10}\right)\frac{U^3}{\ell_0}.$$ 
 \end{Theoreme} 
The main  interest of this result is the fact that we are able to obtain a Kolmogorov type estimate even in the framework of the (damped) Stokes model. The Stokes equations are a  merely a non turbulent model precisely due to the lack of the nonlinear transport term and this raises the problem of the pertinence of this theory for fluids defined in the whole space $\Rt$.
\section{The fractional Navier-Stokes equations}\label{Secc_FracNS}
We study now the fractional damped Navier-Stokes equations introduced in (\ref{Equation_FracNS_Intro}). Our first result establishes the existence of weak solutions and an energy inequality.
\begin{Theoreme}\label{Th_Existence_Fractional_NS} Consider $\vu_0 \in L^{2}(\Rt)$ be a divergence-free initial data and for some $0 <\alpha<4$ let $\vf \in L^2(\Rt)\cap \dot{H}^{-\frac{\alpha}{2}}(\Rt)$ be a divergence free, time independent, external force.  Then, for $0<\beta$ there exists a function
$$\vu=\vu_{\alpha, \beta}\in L^{\infty}_{loc}([0,+\infty[,L^2(\Rt)) \cap L^{2}_{loc}([0,+\infty[, \dot{H}^{\frac{\alpha}{2}} (\Rt)),$$
which is a weak solution of the equation (\ref{Equation_FracNS_Intro}) that satisfies the energy inequality: 
\begin{equation}\label{Energy_inequality_Frac}
\begin{split}
\| \vu(t,\cdot)\|^{2}_{L^2}  \leq \,\, \| \vu_0 \|^{2}_{L^2} -  \, 2 \nu \, \int_{0}^{t}  \left\|\vu(s,\cdot)\right\|^{2}_{\dot{H}^{\frac{\alpha}{2}}} ds + 2 \int_{0}^{t} \int_{\Rt} \vf\cdot\vu(s,\cdot) dx \, ds -2 \beta \int_{0}^{t}\| \vu(s,\cdot)\|^{2}_{L^2} ds.
\end{split}
\end{equation} 
\end{Theoreme}
The proof of this theorem is rather standard and it is based in the hyperviscosity method of Beir\~ao da Vega \cite{Beirao} which introduces the term $-\epsilon \Delta^2$ (this fact explains the technical constraint $0<\alpha<4$). We thus postpone its proof to the Appendix \ref{Secc_AppendixB}. We also refer to \cite{Jarrin1} for another approach to prove this theorem.
\begin{Remarque}\label{Remark_uniqueness_fractional_NS}
Note that, in the undamped Navier-Stokes equation, when $\frac{5}{2} \leq \alpha$ we have the uniqueness of the solution and an energy equality, see \cite{Cholewa}. Moreover, this uniqueness  properties also holds true for the damped equation (\ref{Equation_FracNS_Intro}) since we have a linear damping term.  
\end{Remarque} 
Once we have this existence result, we can see the effect of the damping term $-\beta \vu$ in the definition of the quantity $U$ and by following the same estimates in the proof of Proposition \ref{Prop1} we have:
\begin{Proposition}\label{Prop2}  Within the framework of Theorem \ref{Th_Existence_Fractional_NS} the solutions of the equation (\ref{Equation_FracNS_Intro}) verify: 
\begin{equation*}
\limsup_{T \to +\infty} \frac{1}{T} \int_{0}^{T} \| \vu(t,\cdot) \|^{2}_{L^2} dt \leq \frac{1}{\nu\beta}\,\|\vf\|^{2}_{\dot{H}^{-\frac{\alpha}{2}}}. 
\end{equation*}
\end{Proposition}

\noindent Remark that the previous estimate allows us to define rigorously the quantity $U$ given in (\ref{Def_U_Intro}) in the framework of the fractional and damped Navier-Stokes equation (\ref{Equation_FracNS_Intro}): we can write
$$U= \left( \limsup_{T\to +\infty} \frac{1}{T} \int_{0}^{T} \| \vu(t,\cdot) \|^{2}_{L^2} \frac{dt}{\ell^{3}_{0}}\right)^{\frac{1}{2}}<+\infty.$$
Thus, assuming $U>1$, we obtain the control
\begin{equation}\label{Estimation_UF_frac}
1<U\leq \frac{1}{\sqrt{\nu\beta}}  \frac{\|\vf\|_{\dot{H}^{-\frac{\alpha}{2}}}}{\ell^{\frac{3}{2}}_{0}}.
\end{equation}
Now, for $0<\alpha < 4$ and for $0<\nu$, we define the energy dissipation rate $\mathcal{E}_\alpha$ by the expression
\begin{equation}\label{Def_Epsilon_fractional} 
\mathcal{E}_\alpha= \nu \,    \limsup_{T\to +\infty} \frac{1}{T} \int_{0}^{T} \left\|(-\Delta)^{\frac{\alpha}{4}} \vu(t,\cdot) \right\|^{2}_{L^2} \frac{dt}{\ell^{3}_{0}}.
\end{equation}
Note that by following the same arguments given in the Appendix \ref{Secc_AppendixA} this quantity is well defined. 

\medskip 

With the quantities $U$ and $\mathcal{E}_\alpha$ at our disposal, we will deduce now (formally) a fractional version of the Kolmogorov dissipation law. Indeed, inspired in the relationship (\ref{Kolmogorov_Law_Intro}) we should have 
\begin{equation}\label{Deduccion_Ealpha0}
\mathcal{E}_\alpha \sim \frac{U^{\beta_1}}{\ell^{\beta_2}_{0}},
\end{equation}
where $\beta_1$ and $\beta_2$ are parameters to be fixed and for this we will use a dimensional analysis. First, following some ideas of \cite{Vinegron} (see the formula $(29 b)$ in page $12$) we will assume (only for this analysis) that the velocity $\vu$ is localized at the frecuencies $\frac{1}{\ell_0} < \vert \xi \vert  < \frac{2}{\ell_0}$. Then, by the Bernstein inequalities (see \cite[Chapter $2$]{Grafakos}) and by the expression (\ref{Def_Epsilon_fractional}) we have
\begin{equation}\label{Deduccion_Ealpha1}
\mathcal{E}_{\alpha} \sim \frac{\nu}{\ell^{\alpha}_{0}} U^2.
\end{equation}
Thus, with (\ref{Deduccion_Ealpha0}) and (\ref{Deduccion_Ealpha1}) we can write $\frac{\nu}{\ell^{\alpha}_{0}} U^2 \sim \frac{U^{\beta_1}}{\ell^{\beta_2}_{0}}$, which lead us to the relationship 
$$\nu \sim \frac{U^{\beta_1 -2}}{\ell^{\beta_2-\alpha}_{0}}.$$ 
Now, the term on the right-hand side must have the same physical dimension of the viscosity parameter $\nu$ which is $\frac{length^2}{time}$. Thus, as the characteristic velocity $U$ has a physical dimension of $\frac{length}{time}$ and as $\ell_0$ has a physical dimension of $length$, a simple dimensional analysis shows that we must have $\beta_1 =3$ and $\beta_2 = \alpha-1$. In this way, the Kolmogorov dissipation law in the fractional case writes down as follows: 
\begin{equation}\label{Kolmogorov_Law_fractional}
 \mathcal{E}_\alpha\sim \frac{U^3}{\ell^{\alpha-1}_{0}}, \qquad (Gr\gg 1),
\end{equation} 
Observe that in the case when $\alpha=2$ we obtain the classical damped Navier-Stokes equations (\ref{Eq_damped_NS}), moreover by the expression (\ref{Def_Epsilon_Intro}) and (\ref{Def_Epsilon_fractional}) we have $\mathcal{E} = \mathcal{E}_2$ and we recover the classical Kolmogorov dissipation law (\ref{Kolmogorov_Law_Intro}). \\

The aim of our next result is to prove the relationship (\ref{Kolmogorov_Law_fractional}) above in the setting of the equation (\ref{Equation_FracNS_Intro}). From now on,  only for technical reasons   we shall   consider the range 
$$\frac{3}{7} < \alpha < 3.$$
We will explain below these constraints of the parameter $\alpha$.  Thereafter, by  following the previous analysis done in Section \ref{Secc_Main_theorem} we start by  assuming  that $\vf \in L^{2}(\Rt)\cap \dot{H}^{-\frac{\alpha}{2}}(\Rt)$ satisfies
\begin{equation}\label{Loc_frec_force_Frac}
supp( \widehat{\vf}) \subset \left\{|\xi| < \frac{\mathfrak{c}}{\ell_0}\right\},
\end{equation} 
for some dimensionless constant $0<\mathfrak{c}$, and moreover, for $\ell_0>1$ and $0<\nu <2$, we also assume that  
\begin{equation}\label{Def_F1_fractional_1}
\begin{cases}\vspace{3mm}
\Vert \vf \Vert_{L^2} \gg \ell^{2-\frac{\alpha}{2}}_{0} \qquad \mbox{and} \qquad  \ell_0^{\frac{9}{8}}\|\vf\|_{\dot{H}^{-\frac{\alpha}{2}}}\gg \sqrt{\nu}\|\vf\|_{L^2}^{\frac{7}{4}}, & \quad \mbox{when}\,\, \frac{3}{7}<\alpha < 1, \\ 
\Vert \vf \Vert_{L^2} \gg \ell^{\alpha + \frac{1}{2}}_{0} \qquad \mbox{and} \qquad  \ell_0^{\frac{9}{8}}\|\vf\|_{\dot{H}^{-\frac{\alpha}{2}}}\gg \sqrt{\nu}\|\vf\|_{L^2}^{\frac{7}{4}}, & \quad \mbox{when}\,\, 1\leq \alpha \leq 3. 
\end{cases}
\end{equation} 
This set of conditions will allow us to justify our computations, precisely, the ones given in the proof of Lemma  \ref{Lem_EstimatesFUF_fractional} below where we need to consider the cases $\alpha <1$ and $1 \leq \alpha$ separately. Moreover, these conditions can be simultaneously satisfied by considering, for   example,  a function $\vf$ which is constant in the Fourier level over the ball $B(0, \frac{\mathfrak{c}}{\ell_0})$ and such that 
\begin{equation*}
\begin{cases}\vspace{2mm}
1 \ll \ell^{\frac{7}{2}-\frac{\alpha}{2}}_{0} \leq \Vert \widehat{\vf} \Vert_{L^\infty} \leq \ell^{3+ \frac{2\alpha}{3}}_{0}, & \quad \mbox{when}\,\, \frac{3}{7}\leq \alpha <1, \\
1 \ll \ell^{\alpha+2}_{0} \leq \Vert \widehat{\vf} \Vert_{L^\infty} \leq \ell^{3+ \frac{2\alpha}{3}}_{0}, & \quad \mbox{when}\,\, 1\leq \alpha < 3.  
\end{cases}
\end{equation*}
\begin{Remarque}\label{Rem_LowerUpperBoundAlpha}
Since we assume $\ell_0>1$, the  lower constraint $\frac{3}{7} \leq \alpha$ is needed for the first estimates, while the upper constraint $\alpha < 3$ is crucial for the second estimates. 
\end{Remarque}
We also remark that  (as $\ell_0>1$)  the quantity $\Vert \widehat{\vf} \Vert_{L^\infty}$ is larger as long as $\alpha \to 3^{-}$, \emph{i.e.}, we have a stronger external force.   Finally, let us mention we think that the technical  constraint  $\frac{3}{7} \leq \alpha < 3$ could be removed with another particular (well-prepared) function $\vf$. \\

Once the force is fixed, following (\ref{Def_F}) we define $F= \frac{\| \vf \|_{L^2}}{\ell^{\frac{3}{2}}_{0}}$ and we fix the damping parameter as in (\ref{Def_Beta}), \emph{i.e.} $\beta=F^{\frac{3}{2}}$. In this particular setting we have:
\begin{Theoreme}\label{Th_Fractional_Kolmogorow_Law}  
Let $\vu_0 \in L^2(\Rt)$ be a divergence-free initial data; and consider a length scale $1 \ll \ell_0$ and $0<\nu <2$. Consider also $\vf \in L^2(\Rt)\cap\dot{H}^{-\frac{\alpha}{2}}(\Rt)$ a divergence-free, time independent external force which verifies (\ref{Loc_frec_force_Frac}) and (\ref{Def_F1_fractional_1}).  Let $\vu \in L^{\infty}_{t}L^{2}_{x}\cap (L^{2}_{t})_{loc}\dot{H}^{\frac{\alpha}{2}}_{x}$ be a solution of the equation (\ref{Equation_FracNS_Intro})  associated to $\vu_0, \vf, \alpha$ and $\beta$ (where $\beta$ is defined as in (\ref{Def_Beta})).\\

\noindent If the Grashof number defined in (\ref{Def_Gr}) is large enough (\emph{i.e.} $Gr\gg 1$) then for $\frac{3}{7}\leq \alpha <3$ we have the relationship
\begin{equation}\label{Estimations_TheoremePrincipal_fractional}
\frac{U^3}{\ell_{0}^{\alpha-1}}\sim \mathcal{E}_\alpha, 
\end{equation}
\end{Theoreme}
As pointed out in Remark \ref{Remark_uniqueness_fractional_NS}, in the range $\frac{5}{2}\leq \alpha <3$ we have the uniqueness of weak global solutions of the equation  (\ref{Equation_FracNS_Intro}), but in our current study of the estimates (\ref{Estimations_TheoremePrincipal_fractional}) this property does not play any substantial role.  Nevertheless, we think a deeper study of possible connections between the  Kolmogorov type estimates and the uniqueness (and regularity) of weak solutions of the fractional Navier-Stokes equations could be an interesting matter of further research.\\ 

\noindent {\bf Proof of Theorem \ref{Th_Fractional_Kolmogorow_Law}.} First, we establish an analogous estimate to the one given in Lemma  \ref{Lem_EstimatesFUF}.
\begin{Lemme}\label{Lem_EstimatesFUF_fractional}
Under the hypothesis of the Theorem \ref{Th_Fractional_Kolmogorow_Law} we have the controls
\begin{equation}\label{Estimation_FUF0_fractional}
1\ll F\leq U\leq F^{\frac{3}{2}} \ell_0^{\alpha-1}.
\end{equation}
\end{Lemme}
{\bf Proof.}  By the energy inequality  (\ref{Energy_inequality_Frac}), we also have the  estimate (\ref{Equation_GeneriqueUFU}). From this estimate, by setting $\gamma(2-\gamma)= \frac{1}{U^3}<1$ and as we have $\beta^2=F^3$, we get $F \leq U$.  At this point, we emphasize that by this  last inequality,  by definition of the quantity $F$ (see (\ref{Def_F})) and by  the estimate (\ref{Estimation_UF_frac}), we obtain  the estimate $\frac{\|\vf\|_{L^2}}{\ell_0^{\frac{3}{2}}}\leq \frac{1}{\sqrt{\nu \beta}}\frac{\|\vf\|_{\dot{H}^{-\frac{\alpha}{2}}}}{\ell_0^{\frac{3}{2}}}$ which holds true as long as the second condition in (\ref{Def_F1_fractional_1})  is satisfied with $0<\nu <2$.  

\medskip

On the other hand, by setting now $\gamma(2-\gamma)= \frac{1}{\ell^{2(\alpha-1)}_{0} F^4}<1$ we obtain the estimate $U \leq F^{\frac{3}{2}} \ell^{\alpha-1}_{0}$. Here, it is also worth mentioning  the condition $\frac{1}{\ell^{2(\alpha-1)}_{0} F^4}<1$ is equivalent to $\ell^{2- \frac{\alpha}{2}}_{0} \leq \Vert \vf \Vert_{L^2}$.   When $\frac{3}{7}\leq \alpha <1$, this estimate is precisely  the first condition  in (\ref{Def_F1_fractional_1}). On the other hand, when $1\leq \alpha < 3$, since $\ell_0>1$ by the second condition in (\ref{Def_F1_fractional_1}) we have $\Vert \vf \Vert_{L^2} \gg \ell^{\alpha+\frac{1}{2}}_{0} \geq \ell^{2-\frac{\alpha}{2}}_{0}$.  \hfill $\blacksquare$\\
\begin{itemize}
\item {\bf Lower estimate}. By the energy inequality (\ref{Energy_inequality_Frac}) and by following the same estimates performed to obtain the inequality (\ref{Estim01}), we can write 
$$ \left(\beta-1\right)U^2\leq  \mathcal{E}_\alpha. $$
Then, we have
$$ \left( \frac{\beta \ell^{\alpha-1}_{0}}{U} - \frac{\ell^{\alpha-1}_{0}}{U} \right) \frac{U^3}{\ell^{\alpha-1}_{0}} \leq \mathcal{E}_{\alpha},$$ 
where we must find a lower bound for each expression in the term  on the left-hand side. As $\beta^2=F^3$, for the first term we have $ \frac{\beta \ell^{\alpha-1}_{0}}{U}  =  \frac{F^{\frac{3}{2}} \ell^{\alpha-1}_{0}}{U}$. Moreover, by the estimate (\ref{Estimation_FUF0_fractional}) we have $1 \leq \frac{F^{\frac{3}{2}} \ell^{\alpha-1}_{0}}{U}$. On the other hand, for the second term, as by the estimate  (\ref{Estimation_FUF0_fractional})  we have $F \leq U$ then we obtain  $- \frac{\ell^{\alpha-1}_{0}}{F} \leq - \frac{\ell^{\alpha-1}_{0}}{U}$. Then, as $Gr= \frac{F \ell^{3}_{0}}{\nu^2}$ we can write $- \frac{\ell^{\alpha-1}_{0}}{F} = -\frac{\ell^{\alpha+2}_{0}}{\nu^2\, Gr}$. 

\medskip

We thus have the lower bound $\left( 1 -  \frac{\ell^{\alpha+2}_{0}}{\nu^2\, Gr}\right) \leq \left( \frac{\beta \ell^{\alpha-1}_{0}}{U} - \frac{\ell^{\alpha-1}_{0}}{U} \right)$ from which we get the lower estimate
$$ \left( 1 -  \frac{\ell^{\alpha+2}_{0}}{\nu^2\, Gr}\right)  \frac{U^3}{\ell^{\alpha-1}_{0}}\leq \mathcal{E}_\alpha.$$
Now, we shall need the following estimate
\begin{equation}\label{Estim_Gr_ell_nu_fractional}
\frac{\ell^{\alpha+2}_{0}}{\nu^2\, Gr} \leq \frac{1}{2},
\end{equation}
which follows from the condition (\ref{Def_F1_fractional_1})  according to the different cases of the parameter $\alpha$. Indeed, since $Gr= \frac{F \ell^{3}_{0}}{\nu^2}$ the estimate above is equivalent to the estimate $2 \ell^{\alpha+\frac{1}{2}}_{0} \leq \Vert \vf \Vert_{L^2}$. When $\frac{3}{7}\leq \alpha <1$ this estimate is satisfied by the first condition in (\ref{Def_F1_fractional_1}) since we have $\ell_0>1$ and $\Vert \vf \Vert_{L^2}\gg \ell^{2-\frac{\alpha}{2}} \gg \ell^{\alpha+\frac{1}{2}}_{0}$. Thereafter, when $1\leq \alpha < 3$ this estimate is precisely given by the first  condition in (\ref{Def_F1_fractional_1}). 

\medskip

With the estimate (\ref{Estim_Gr_ell_nu_fractional}) at our disposal,  we obtain the lower estimate 
$$ \frac{1}{2} \frac{U^3}{\ell^{\alpha-1}_{0}} \leq \mathcal{E}_\alpha. $$
\item {\bf Upper estimate.} Always by the energy inequality (\ref{Energy_inequality_Frac}) and by following the same estimates performed to obtain the inequality (\ref{Estim02}),  for $0<\rho$ we have
$$ 2\mathcal{E}_\alpha \leq  \frac{1}{\rho}F^2+\rho U^2. $$
Then,  we set $\rho = \frac{\ell^{\alpha-1}_{0}}{10\, U}$ to obtain 
$$ \mathcal{E}_\alpha  \leq \frac{10 \, U}{\ell^{\alpha-1}_{0}} F^2 +  \frac{\ell^{\alpha-1}_{0} }{10} U. $$
By the lower estimate in (\ref{Estimation_FUF0_fractional}) we have $F^2 \leq U^2$ and then we get
$$ \mathcal{E}_\alpha  \leq  10 \frac{U^3}{\ell^{\alpha-1}_{0}} +  \frac{\ell^{\alpha-1}_{0} }{10} U = \frac{U^3}{\ell^{\alpha-1}_{0}} \left( 10 + \frac{\ell^{2(\alpha-1)}_{0}}{10 U^2}  \right)  \leq \frac{U^3}{\ell^{\alpha-1}_{0}} \left( 10 + \frac{\ell^{2(\alpha-1)}_{0}}{10 F^2}  \right).$$ 
Finally, always by the estimate (\ref{Estim_Gr_ell_nu_fractional})  we have $\frac{\ell^{2(\alpha-1)}_{0}}{ F^2}= \frac{\ell^{2\alpha+4}_{0}}{\nu^4 Gr^2} \leq \frac{1}{4}$, which yields the upper estimate
$$ \mathcal{E}_\alpha \leq  \frac{401}{40}\, \frac{U^3}{\ell^{\alpha-1}_{0}}.$$     
\end{itemize} 
Theorem \ref{Th_Fractional_Kolmogorow_Law} is now proven.  \hfill $\blacksquare$
\appendix
\section{The energy dissipation rate}\label{Secc_AppendixA}
We prove here that que quantity $\mathcal{E}$ given in (\ref{Def_Epsilon_Intro}) is well defined. Indeed, by the energy inequality (\ref{Inegalite_Energie_Intro}) we just write 
\begin{equation*}
\begin{split}
2 \nu \int_{0}^{T} \Vert \vu(t,\cdot)\Vert^{2}_{\dot{H}^1} dt \leq &\, \Vert \vu \Vert^{2}_{L^2} + 2 \int_{0}^{T}\langle \vf, \vu(t,\cdot)\rangle_{\dot{H}^{-1} \times \dot{H}^1}  dt   \leq \Vert \vu \Vert^{2}_{L^2}  + 2\int_{0}^{T} \Vert \vf \Vert_{\dot{H}^{-1}} \Vert \vu(t,\cdot)\Vert_{\dot{H}^1} dt \\
\leq &\, \Vert \vu \Vert^{2}_{L^2} + \frac{T}{\nu} \Vert \vf \Vert^{2}_{\dot{H}^{-1}} + \nu \int_{0}^{T}\Vert \vu(t,\cdot)\Vert^{2}_{\dot{H}^1} dt,
\end{split}
\end{equation*}
hence we have
$$ \nu \int_{0}^{T} \Vert \vu(t,\cdot)\Vert^{2}_{\dot{H}^1} dt \leq \Vert \vu \Vert^{2}_{L^2} + \frac{T}{\nu} \Vert \vf \Vert^{2}_{\dot{H}^{-1}}.$$
Then, we divide each term by $T$ and by $\ell^{3}_{0}$, and moreover, by letting $T \to +\infty$ we have  the control $\mathcal{E} \leq \frac{\Vert \vf \Vert^{2}_{\dot{H}^{-1}}}{\nu\, \ell^{3}_{0}} < +\infty$. 
\section{Existence of weak solutions for the system (\ref{Equation_FracNS_Intro}) }\label{Secc_AppendixB}
There are several methods to obtain a weak solution of the equation 
\begin{equation*}
\left\lbrace \begin{array}{lc}\vspace{3mm}
\partial_t\vu=- \nu (-\Delta)^{\frac{\alpha}{2}} \vu-\P((\vu\cdot \vec{\nabla}) \vu)+\vf-\beta \vu, \quad div(\vu)=0,\quad \alpha, \beta, \, \nu>0,  &\\
\vu(0,\cdot)=\vu_0.&\\
\end{array}\right.
\end{equation*}
Perhaps the simpler way is to follow the hyperviscosity approach given by \cite{Beirao} and to consider the following model
\begin{equation}\label{Equation_FracNS_hyperviscosity}
\partial_t\vu=-\epsilon\Delta^2\vu- \nu (-\Delta)^{\frac{\alpha}{2}} \vu-\P((\vu\cdot \vec{\nabla}) \vu)+\vf-\beta \vu, \quad \varepsilon>0.
\end{equation}
In the sequel, we shall follow the program  given in \cite{Beirao} to construct global in time weak solutions of the equation (\ref{Equation_FracNS_Intro}).  We shall mainly  detail the  estimates involving the terms $-\beta \vu$ and $-\nu (-\Delta)^{\frac{\alpha}{2}}\vu$ which are different from the classical case. \\

\noindent {\bf First step: Weak solutions for the equation (\ref{Equation_FracNS_hyperviscosity})}.  We recall that the operator $e^{-\varepsilon t \Delta^2}$ is defined in the Fourier level as $ e^{-\varepsilon t \vert \xi \vert^4} \widehat{\varphi}$ for all $\varphi$ belonging to the Schwartz class.  For $T>0$ small enough, we will solve the following fixed point problem
\begin{equation}\label{Fixed-point-FracNS_hyperviscosity}
\begin{split}
\vu(t,x)= &\,  e^{-\varepsilon t \Delta^2} \vu_0(x) + \int_{0}^{t} e^{-\varepsilon (t-s) \Delta^2} \vf(x) ds  -  \int_{0}^{t} e^{-\varepsilon (t-s) \Delta^2} ((\vu \cdot \vn) \vu)(s,x) ds \\
&\, - \int_{0}^{t} e^{-\varepsilon (t-s) \Delta^2} (\beta I_d+\nu (-\Delta)^{\frac{\alpha}{2}})\vu(s,x) ds 
\end{split}
\end{equation}
in the space $L^{\infty}([0,T], L^2(\Rt))\cap L^{2}([0,T], \dot{H}^2(\Rt))$ with the norm $\Vert \cdot \Vert_{T}= \Vert \cdot \Vert_{L^{\infty}_{t}L^{2}_{x}} + \sqrt{\varepsilon}\Vert \cdot \Vert_{L^{2}_{t} \dot{H}^{2}_{x}}$.  We thus write 
\begin{equation}
\begin{split}
\Vert \vu  \Vert_T= &\, \underbrace{\left\Vert e^{-\varepsilon t \Delta^2} \vu_0 \right\Vert_T}_{(A)} + \underbrace{\left\Vert \int_{0}^{t} e^{-\varepsilon (t-s) \Delta^2} \vf(\cdot) ds \right\Vert_{T}}_{(B)}  + \underbrace{ \left\Vert   \int_{0}^{t} e^{-\varepsilon (t-s) \Delta^2} \P ((\vu \cdot \vn) \vu)(s,\cdot) ds \right\Vert_{T}}_{(C)} \\
&\, + \underbrace{ \left\Vert  \int_{0}^{t} e^{-\varepsilon (t-s) \Delta^2} (\beta I_d+\nu (-\Delta)^{\frac{\alpha}{2}})\vu(s,\cdot) ds \right\Vert_{T}}_{(D)}.  
\end{split}
\end{equation}
For the term $(A)$, by \cite{PGL}, Theorem $18.2$  (page $352$),  we have the estimate $\left\Vert e^{-\varepsilon t \Delta^2} \vu_0 \right\Vert_T \leq C\Vert \vu_0 \Vert_{L^2}$. In order to study the  term $(B)$, we  can adapt the estimates given in  \cite{PGL}, Theorem $18.2$, to obtain  the following: 
\begin{Lemme} Let $\vg \in L^{2}_{t}L^{2}_{x}$. Then we have 
\begin{equation}\label{Estim-g-0}
\left\Vert \int_{0}^{t} e^{-\varepsilon (t-s) \Delta^2} \vg(s,\cdot) ds \right\Vert_{T} \leq C_1 \,  \Vert \vg \Vert_{L^{\infty}_{t} L^{2}_{x}}, \quad \text{with}\quad C_1=  C\left(  T^{\frac{1}{2}} + \sqrt{\varepsilon} T^{\frac{1}{2}} + \frac{1}{\sqrt{\varepsilon}}\right).
\end{equation} 
\end{Lemme}	
{\bf Proof.}  For the first term of the norm $\Vert \cdot \Vert_{T}$, for $0<t \leq T$ we have 
\begin{equation}\label{Estim-g-1}
\left\Vert \int_{0}^{t} e^{-\varepsilon (t-s) \Delta^2} \vg(\cdot) ds \right\Vert_{L^2} \leq  \int_{0}^{t} \Vert e^{-\varepsilon (t-s) \Delta^2} \vg(\cdot)  \Vert_{L^2} ds \leq \int_{0}^{t} \Vert \vg(s,\cdot) \Vert_{L^2} ds\leq   T^{\frac{1}{2}}\, \Vert \vg \Vert_{L^{2}_{t}L^{2}_{x}}. 
\end{equation}
On the other hand, for the second term of the norm $\Vert \cdot \Vert_{T}$, we have
\begin{equation*}
\begin{split}
\sqrt{\varepsilon}\, \left\Vert  \int_{0}^{t} e^{-\varepsilon (t-s) \Delta^2} \vg(\cdot) ds \right\Vert_{L^{2}_{t}\dot{H}^{2}_{x}}  \leq &\,   \sqrt{\varepsilon}\, \left\Vert  \int_{0}^{t} e^{-\varepsilon (t-s) \Delta^2} \vg(\cdot) ds  \right\Vert_{L^{2}_{t}H^{2}_{x}} \\
 \leq &\, \sqrt{\varepsilon}\, \left\Vert  (I_d -\Delta) \left(  \int_{0}^{t} e^{-\varepsilon (t-s) \Delta^2} \vg(\cdot) ds  \right) \right\Vert_{L^{2}_{t} L^{2}_{x}}\\
  \leq &\, \sqrt{\varepsilon}\, \left\Vert  (I_d -\Delta)^2 \left(  \int_{0}^{t} e^{-\varepsilon (t-s) \Delta^2} (I_d -\Delta)^{-1} \vg(\cdot) ds  \right) \right\Vert_{L^{2}_{t} L^{2}_{x}}\\
 \leq &\,  C \sqrt{\varepsilon}\, \left\Vert   \left(  \int_{0}^{t} e^{-\varepsilon (t-s) \Delta^2} (I_d -\Delta)^{-1} \vg(\cdot) ds  \right) \right\Vert_{L^{2}_{t} L^{2}_{x}} \\
 &\,+ C \sqrt{\varepsilon}\, \left\Vert  \Delta^2 \left(  \int_{0}^{t} e^{-\varepsilon (t-s) \Delta^2} (I_d -\Delta)^{-1} \vg(\cdot) ds  \right) \right\Vert_{L^{2}_{t} L^{2}_{x}}.
\end{split}
\end{equation*}
For the first term on the right-hand side, by (\ref{Estim-g-1}) we can write  
\begin{equation*}
C \sqrt{\varepsilon}\, \left\Vert   \left(  \int_{0}^{t} e^{-\varepsilon (t-s) \Delta^2} (I_d -\Delta)^{-1} \vg(\cdot) ds  \right) \right\Vert_{L^{2}_{t} L^{2}_{x}} \leq C\sqrt{\varepsilon} T^{\frac{1}{2}}  \Vert (I_d -\Delta)^{-1}  \vg \Vert_{L^{2}_{t} L^{2}_{x}} \leq  C\sqrt{\varepsilon} T^{\frac{1}{2}} \Vert  \vg \Vert_{L^{2}_{t} L^{2}_{x}}
\end{equation*}
Moreover, for the second term on the right-hand side, by the maximal regularity of the bi-Laplacian operator  we can write $ C \sqrt{\varepsilon}\, \left\Vert  \Delta^2 \left(  \int_{0}^{t} e^{-\varepsilon (t-s) \Delta^2} (I_d -\Delta)^{-1} \vg(\cdot) ds  \right) \right\Vert_{L^{2}_{t} L^{2}_{x}} \leq \frac{C}{\sqrt{\varepsilon}} \Vert \vg \Vert_{L^{2}_{t}L^{2}_{x}}$.  By gathering all these estimates we get the wished estimate (\ref{Estim-g-0}). \finpv 
\medskip

Once we have the estimate (\ref{Estim-g-0}) at our disposal, be setting $\vg = \vf$ (recall that $\vf \in L^2(\Rt)$ is time-independent) for the term $(B)$ we have  $\left\Vert \int_{0}^{t} e^{-\varepsilon (t-s) \Delta^2} \vf(\cdot) ds \right\Vert_{T} \leq C_1 \, T^{\frac{1}{2}}   \Vert \vf \Vert_{L^{2}}$. We thus can write
\begin{equation}\label{Estim-data}
\left\Vert e^{-\varepsilon t \Delta^2} \vu_0  + \int_{0}^{t} e^{-\varepsilon (t-s) \Delta^2} \vf(\cdot) ds  \right\Vert_T \leq C \Vert \vu_0 \Vert_{L^2} + C_1 T^{\frac{1}{2}} \Vert \vf \Vert_{L^2}. 
\end{equation}
In order to study the  term $(C)$, for a constant $C_B=  C \varepsilon^{-\frac{3}{8}}\left( \left( \frac{T}{\varepsilon}\right)^{\frac{1}{4}}+\sqrt{\varepsilon} T^{\frac{3}{4}}+\frac{1}{\sqrt{\varepsilon}}\right)\left( 1 +T^{\frac{1}{2}}\right)$ by \cite{PGL}, Theorem $18.2$,  we have 
\begin{equation}\label{Estim-bilinear}
 \left\Vert   \int_{0}^{t} e^{-\varepsilon (t-s) \Delta^2} \P ((\vu \cdot \vn) \vu)(s,\cdot) ds \right\Vert_{T} \leq C_B T^{\frac{1}{8}} \Vert \vu \Vert^{2}_{T}.
\end{equation}
Finally, to study  the term $(D)$ we write
\begin{equation*}
\begin{split}
 \left\Vert  \int_{0}^{t} e^{-\varepsilon (t-s) \Delta^2} (\beta I_d+\nu (-\Delta)^{\frac{\alpha}{2}})\vu(s,\cdot) ds \right\Vert_{T} \leq  &\, \beta  \left\Vert  \int_{0}^{t} e^{-\varepsilon (t-s) \Delta^2} \vu(s,\cdot) ds \right\Vert_{T}  \\
 &\, + \nu   \left\Vert  \int_{0}^{t} e^{-\varepsilon (t-s) \Delta^2}(-\Delta)^{\frac{\alpha}{2}}\vu(s,\cdot) ds \right\Vert_{T}
 \end{split}
\end{equation*}
For the first term on the right-hand side, we use again the estimate (\ref{Estim-g-0}) (with $\vg=  \vu$) to write 
\begin{equation}
\beta \left\Vert  \int_{0}^{t} e^{-\varepsilon (t-s) \Delta^2} \vu(s,\cdot) ds \right\Vert_{T} \leq \beta C_1  \Vert \vu \Vert_{L^{2}_{t}L^{2}_{x}} \leq \beta C_1 T^{\frac{1}{2}} \Vert \vu \Vert_{L^{\infty}_{t}L^{2}_{x}} \leq \beta C_1 T^{\frac{1}{2}} \Vert \vu \Vert_{T}.  
\end{equation}
For the second term on the right-hand side, we shall need the following:
\begin{Lemme}\label{Estim-Linear} Let $0<\alpha <4$ and let $\vu \in L^{\infty}_{t}L^{2}_{x} \cap L^{2}_{t}\dot{H}^{2}_{x}$. The following estimates hold:
\begin{enumerate}
\item[1)] When $0<\alpha \leq 2$ we have $\ds{ \left\Vert  \int_{0}^{t} e^{-\varepsilon (t-s) \Delta^2}(-\Delta)^{\frac{\alpha}{2}}\vu(s,\cdot) ds \right\Vert_{T} \leq C\left( \left( \frac{T}{\varepsilon}\right)^{1-\frac{\alpha}{4}}+ \sqrt{\varepsilon} T^{\frac{1}{2}} \right) \Vert \vu \Vert_{T}}$.
\item[2)] When $2<\alpha <4$ we have $\ds{ \left\Vert  \int_{0}^{t} e^{-\varepsilon (t-s) \Delta^2}(-\Delta)^{\frac{\alpha}{2}}\vu(s,\cdot) ds \right\Vert_{T} C\left( \left( \frac{T}{\varepsilon}\right)^{1-\frac{\alpha}{4}}+ \varepsilon^{\frac{3}{4} -\frac{\alpha}{8}}\,T^{\frac{6-\alpha}{8}} \right) \Vert \vu \Vert_{T} }$.
\end{enumerate}		
\end{Lemme}	
{\bf Proof.} For $0<\alpha <4$ we estimate   the first term  in the norm $\Vert \cdot \Vert_{T}$. By duality and by the Cauchy-Schwarz inequality we have 
\begin{equation*}
\begin{split}
 \left\Vert  \int_{0}^{t} e^{-\varepsilon (t-s) \Delta^2}(-\Delta)^{\frac{\alpha}{2}}\vu(s,\cdot) ds \right\Vert_{L^2} =&\, \sup_{\Vert \vphi \Vert_{L^2} \leq 1 } \left\vert  \int_{\Rt} \left( \int_{0}^{t} e^{-\varepsilon(t-s)\Delta^2} (-\Delta)^{\frac{\alpha}{2}} \vu(s,x) ds\right) \cdot \vphi(x) dx \right\vert  \\
=&\, \sup_{\Vert \vphi \Vert_{L^2} \leq 1 } \left\vert  \int_{0}^{t} \int_{\Rt} \vu(s,x)\cdot e^{-\varepsilon(t-s)\Delta^2} (-\Delta)^{\frac{\alpha}{2}} \vphi(x) dx \, ds  \right\vert\\
\leq &\, \Vert \vu \Vert_{L^{\infty}_{t} L^{2}_{x}} \,  \sup_{\Vert \vphi \Vert_{L^2} \leq 1 }  \, \int_{0}^{t} \Vert e^{-\varepsilon(t-s)\Delta^2} (-\Delta)^{\frac{\alpha}{2}} \vphi \Vert_{L^2}, 
\end{split} 
\end{equation*}
where we must study the integral above. By the Plancherel identity we write
\begin{equation*}
\begin{split}
\int_{0}^{t} \Vert e^{-\varepsilon(t-s)\Delta^2} (-\Delta)^{\frac{\alpha}{2}} \vphi \Vert_{L^2} =&\, \int_{0}^{t} \left( \int_{\Rt} e^{-2  \vert (\varepsilon(t-s))^{\frac{1}{4}} \xi \vert^4} \vert (\varepsilon (t-s))^{\frac{1}{4}} \xi \vert^{2 \alpha} \vert \widehat{\vphi}(\xi) \vert^{2} d \xi \right)^{\frac{1}{2}} (\varepsilon (t-s))^{-\frac{\alpha}{4}} \, ds\\
\leq &\, C \Vert \vphi \Vert_{L^2} \left( \frac{t}{\varepsilon} \right)^{1-\frac{\alpha}{4}}  \leq C \left( \frac{T}{\varepsilon} \right)^{1-\frac{\alpha}{4}}. 
\end{split}
\end{equation*}
We thus obtain 
\begin{equation}\label{Estim-ter-lin-1}
  \left\Vert  \int_{0}^{t} e^{-\varepsilon (t-s) \Delta^2}(-\Delta)^{\frac{\alpha}{2}}\vu(s,\cdot) ds \right\Vert_{L^2} \leq  \, C \left( \frac{T}{\varepsilon} \right)^{1-\frac{\alpha}{4}} \, \Vert \vu \Vert_{T}. 
\end{equation}
For the second term in the norm $\Vert \cdot \Vert_{T}$, always by  the Plancherel identity  we can write
\begin{equation*}
\begin{split}
\sqrt{\varepsilon} \left\Vert  \int_{0}^{t} e^{-\varepsilon (t-s) \Delta^2}(-\Delta)^{\frac{\alpha}{2}}\vu(s,\cdot) ds \right\Vert_{L^{2}_{t}\dot{H}^{2}} = \sqrt{\varepsilon}  \left\Vert  (-\Delta)^{1+\frac{\alpha}{2}} \left( \int_{0}^{t}  e^{-\varepsilon(t-s)\Delta^2} \vu(s,\cdot) ds \right) \right\Vert_{L^{2}_{t}L^{2}_{x}},  
\end{split}
\end{equation*}
where we must treat the cases $0<\alpha \leq 2$ and $2<\alpha <4$ separately. For the first case, by the maximal regularity of the bi-Laplacian operator we have 
\begin{equation}\label{Estim-ter-lin2}
\begin{split}
\sqrt{\varepsilon} \left\Vert  (-\Delta)^{1+\frac{\alpha}{2}} \left( \int_{0}^{t}  e^{-\varepsilon(t-s)\Delta^2} \vu(s,\cdot) ds \right) \right\Vert_{L^{2}_{t}L^{2}_{x}} \leq \sqrt{\varepsilon} C\, \Vert \vu \Vert_{L^{2}_{t}L^{2}_{x}} \leq  \sqrt{\varepsilon} C\, T^{\frac{1}{2}}\,  \Vert \vu \Vert_{L^{\infty}_{t}L^{2}_{x}} \leq  \sqrt{\varepsilon} C\, T^{\frac{1}{2}}\,  \Vert \vu \Vert_{T}. 
\end{split}
\end{equation}
For the second case, as $2<\alpha <4$ we have $1<\frac{\alpha}{2}<2$ and for $0<s<1$ we write $\frac{\alpha}{2}=1+s$. Then, always by the maximal regularity of the bi-Laplacian operator  we can write 
\begin{equation*}
\begin{split}
\sqrt{\varepsilon} \left\Vert  (-\Delta)^{1+\frac{\alpha}{2}} \left( \int_{0}^{t}  e^{-\varepsilon(t-s)\Delta^2} \vu(s,\cdot) ds \right) \right\Vert_{L^{2}_{t}L^{2}_{x}}= &\, \sqrt{\varepsilon} \left\Vert  (-\Delta)^{2} \left( \int_{0}^{t}  e^{-\varepsilon(t-s)\Delta^2} (-\Delta)^{s} \vu(s,\cdot) ds \right) \right\Vert_{L^{2}_{t}L^{2}_{x}}\\
\leq &\, \sqrt{\varepsilon}\, C \Vert (-\Delta)^s \vu \Vert_{L^{2}_{t} L^{2}_{x}}. 
\end{split}
\end{equation*}
Then, since  $\vu \in L^{\infty}_{t} L^{2}_{x} \cap L^{2}_{t}\dot{H}^{2}_{x}$ and since $0<s<1$, for $p= \frac{8}{\alpha-2}>2$ (recall that $2<\alpha<4$) by the interpolation inequalities with $\theta=1-\frac{s}{2}\in (\frac{1}{2}, 1)$, we have
\begin{equation*}
\begin{split}
\sqrt{\varepsilon}C\,\Vert (-\Delta)^s \vu \Vert_{L^{2}_{t} L^{2}_{x}}  \leq &\, \sqrt{\varepsilon}C\,\Vert  \vu \Vert_{L^{2}_{t} \dot{H}^{2s}_{x}} \leq \sqrt{\varepsilon}C T^{\frac{1}{2}-\frac{1}{p}} \, \Vert \vu \Vert_{L^{p}_{t} \dot{H}^{2s}_{x}} \leq \sqrt{\varepsilon}C \, T^{\frac{1}{2}-\frac{1}{p}}\,  \Vert \vu \Vert^{\theta}_{L^{\infty}_{t} L^{2}_{x}}\, \Vert \vu \Vert^{1-\theta}_{L^{2}_{t} \dot{H}^{2}_{x}}  \\
\leq  &\, \sqrt{\varepsilon}C\,  (\sqrt{\varepsilon})^{\theta-1}T^{\frac{1}{2}-\frac{1}{p}}\,  \Vert \vu \Vert^{\theta}_{L^{\infty}_{t} L^{2}_{x}}\, ( \sqrt{\varepsilon}\Vert \vu \Vert_{L^{2}_{t} \dot{H}^{2}_{x}})^{1-\theta} \leq  (\sqrt{\varepsilon})^{\theta}C\, T^{\frac{1}{2}-\frac{1}{p}}\, \Vert \vu \Vert_{T}\\
=&\, \varepsilon^{\frac{3}{4} -\frac{\alpha}{8}} C\, T^{\frac{6-\alpha}{8}}\, \Vert \vu \Vert_{T}. 
\end{split}
\end{equation*}
Then we obtain 
\begin{equation}\label{Estim-ter-lin3}
\sqrt{\varepsilon} \left\Vert  (-\Delta)^{1+\frac{\alpha}{2}} \left( \int_{0}^{t}  e^{-\varepsilon(t-s)\Delta^2} \vu(s,\cdot) ds \right) \right\Vert_{L^{2}_{t}L^{2}_{x}} \leq \varepsilon^{\frac{3}{4} -\frac{\alpha}{8}} C\, T^{\frac{6-\alpha}{8}}\, \Vert \vu \Vert_{T}.
\end{equation}
The first point stated in this lemma follows from the estimates (\ref{Estim-ter-lin-1}) and (\ref{Estim-ter-lin2}), while the second point follows from the estimates (\ref{Estim-ter-lin-1}) and (\ref{Estim-ter-lin3}). \finpv 
\medskip

Once we have the estimates (\ref{Estim-data}), (\ref{Estim-bilinear}) and the estimates given in Lemma \ref{Estim-Linear} at our disposal, for a time $T>0$ and for $\varepsilon>0$ by the Banach contraction principle we obtain $\vu_\varepsilon \in L^{\infty}([0,T],L^{2}(\Rt))\cap L^{2}([0,T], \dot{H}^2(\Rt))$ a local solution of the equations (\ref{Fixed-point-FracNS_hyperviscosity}).\\

This solution also verifies the hyperviscosity equation  (\ref{Equation_FracNS_hyperviscosity}) and by standard computations, for $t\geq 0$ we can obtain the following energy inequality 
\begin{equation}\label{Energy-inequality-hyperviscosity}
\begin{split}
&\, \Vert \vu_{\varepsilon}(t,\cdot)\Vert^{2}_{L^2} + 2\nu \int_{0}^{t}\Vert \vu(s,\cdot)\Vert^{2}_{\dot{H}^{\frac{\alpha}{2}}} ds + 2\varepsilon \int_{0}^{t} \Vert \vu_{\varepsilon}(s,\cdot)\Vert^{2}_{\dot{H}^2} ds + 2\beta \int_{0}^{t}\Vert \vu_{\varepsilon}(s,\cdot)\Vert^{2}_{L^2}ds \\
\leq &\, \Vert \vu_0 \Vert^{2}_{L^2}+ 2\int_{0}^{t} \int_{\Rt} \vf (x)\cdot \vu_{\varepsilon}(s,x) dx ds. 
\end{split}
\end{equation}
We thus deduce the following control 
\begin{equation}\label{Energy-control-hyperviscosity}
\Vert \vu_{\varepsilon}(t,\cdot)\Vert^{2}_{L^2} + \nu \int_{0}^{t}\Vert \vu(s,\cdot)\Vert^{2}_{\dot{H}^{\frac{\alpha}{2}}} ds + 2\varepsilon \int_{0}^{t} \Vert \vu_{\varepsilon}(s,\cdot)\Vert^{2}_{\dot{H}^2} ds 
\leq \, \Vert \vu_0 \Vert^{2}_{L^2}+ t \frac{\Vert \vf \Vert^{2}_{\dot{H}^{-\frac{\alpha}{2}}}}{\nu}\\
\end{equation}
 which  allows us to extend the local solution $\vu_{\varepsilon}$ to the whole interval $[0,+\infty[$. \\[5mm]
{\bf Second step: Passage to the limit and weak solutions for the equation (\ref{Equation_FracNS_Intro})}.  By the control (\ref{Energy-control-hyperviscosity}) the family $(\vu_{\varepsilon})_{\varepsilon>0}$ of solutions of the equation  (\ref{Equation_FracNS_hyperviscosity}) is uniformly bounded (respect to the parameter $\varepsilon$) in the space $(L^{\infty}_{t})_{loc}L^{2}_{x}\cap  (L^{2}_{t})_{loc}\dot{H}^{\frac{\alpha}{2}}_{x}$. Consequently, there exists $\vu \in  (L^{\infty}_{t})_{loc}L^{2}_{x}\cap  (L^{2}_{t})_{loc}\dot{H}^{\frac{\alpha}{2}}_{x}$  and there exits a  subsequence $(\varepsilon_n)_{n \in \mathbb{N}}$ such that for all $0<T<+\infty$ we have $\vu_{\varepsilon_n} \to \vu$ in the weak-*  topology of the spaces $L^{\infty}([0,T],L^{2}(\Rt))$ and $L^{2}([0,T], \dot{H}^{\frac{\alpha}{2}}(\Rt))$.  Moreover, by the  Rellich-Lions lemma (see \cite{PGL}, Theorem $12.1$) we also have that $\vu_{\varepsilon_n} \to \vu$ in the strong topology of the space $L^{2}_{loc}([0,+\infty[ \times \Rt)$. From these convergences we can deduce that the sequence $\left( \P\left( (\vu_{\varepsilon_n} \cdot \vn) \vu_{\varepsilon_n}\right) \right)_{n \in \mathbb{N}}$ converges to $\P\left( (\vu  \cdot \vn) \vu \right)$ in the weak-* topology of the space $(L^{2}_{t})_{loc} H^{\frac{\alpha}{2}-4}_{x}$ and then the limit $\vu$ is a weak solution of the equation (\ref{Eq_damped_NS}).\\ 

Finally, from the energy inequality (\ref{Energy-inequality-hyperviscosity}) by applying  standard methods (see \cite{PGL}, Theorem $18.2$) we obtain the energy inequality (\ref{Energy_inequality_Frac}). Theorem \ref{Th_Existence_Fractional_NS} is proven. \finpv 
\section{A Kolmogorov type estimate for another choice of the damping parameter}\label{Secc_AppendixC}
By following some ideas of \cite{Chamorro} and \cite[Chapter $1$]{Jarrin}, from  the length  $\ell_0$ and the viscosity parameter $\nu$ we set now the damping parameter 
\begin{equation}\label{Def_Beta_2}
\beta = \frac{\nu}{\ell^{2}_{0}}. 
\end{equation}
We remark that this value of the damping parameter is different to the one given in (\ref{Def_Beta}) since it does not depend on the quantity $\Vert \vf \Vert_{L^2}$. Moreover, this choice of the damping parameter can be motivated from the space-periodic setting as follows: when we consider the classical Navier-Stokes equations with space-periodic conditions in the box $\Omega=[0,\ell_0]^3$, for $\vu \in L^{\infty}_{t}L^{2}_{x}\cap L^{2}_{t}\dot{H}^{1}_{x}$ a generic Leray solution,  by the well-known energy inequality and by the \emph{Poincaré's inequality} we have the control:
\[  \limsup_{T\to +\infty} \frac{1}{T} \int_{0}^{T} \Vert \vu(t,\cdot)\Vert^{2}_{L^{2}(\Omega)}\, dt  \leq   \frac{\ell^{2}_{0}}{\nu^2}\Vert \vf \Vert^{2}_{\dot{H}^{-1}}, \]
see \cite[Section $1.2.3$]{Jarrin} for all the details. Thus, when comparing this control to the one given in Proposition \ref{Prop1} (in the setting of the whole space) we observe that  the identity (\ref{Def_Beta_2})  naturally appears as a choice for the damping parameter $\beta$. In this sense, we remark that the damping term $-\beta \vu$ in the equations (\ref{Eq_damped_NS}) acts as compensation for the lack of the Poincaré inequality in the whole space $\Rt$.

\medskip

In the following theorem we show that the  value of the damping parameter  given in the identity  (\ref{Def_Beta_2}) also allows to prove an estimate from above according to  the Kolmogorov dissipation law (\ref{Kolmogorov_FRAC_Law_Intro1}) in the deterministic setting of the damped Navier-Stokes equations (\ref{Eq_damped_NS}).  
\begin{Theoreme}\label{Th_Kolmogorow_Law_Damping2} With the same hypothesis of Theorem \ref{Th_Kolmogorow_Law},  consider now the damping parameter $\beta$ given in (\ref{Def_Beta_2}). Moreover, let  the Grashof number $Gr$ and the Reynolds number $Re$ defined in (\ref{Def_Gr}) and (\ref{Def_Re_Intro}) respectively. \\
	
There exists  a constant $0<{\bf c}$, which do not depend on any physical quantity, such that if $2Gr\leq Re$ then the following estimate holds:
\begin{equation}\label{Estimate_above_Kolmogorov_Damping2}
 \mathcal{E}\leq {\bf c} \left( \frac{1}{Gr} + 1\right)\, \frac{U^3}{\ell_0},
\end{equation}
where the quantities $U$ and $\mathcal{E}$ are given in (\ref{Def_U_Intro}) and (\ref{Def_Epsilon_Intro}) respectively. 
\end{Theoreme}	 

\medskip

We emphasize that here we need to assume the (technical) control $2Gr \leq Re$ which, to our knowledge, does not follow from  the hypothesis of this theorem. In contrast, when we set the damping parameter $\beta$ as is (\ref{Def_Beta}),  this type of control can be obtained from the hypothesis of Theorem \ref{Th_Kolmogorow_Law}. See the Corollary \ref{Coro_Relation_ReGr} above.  \\

\noindent {\bf Proof of Theorem  \ref{Th_Kolmogorow_Law_Damping2}.} The proof of  the  estimate (\ref{Estimate_above_Kolmogorov_Damping2}) is based on the following inequalities:
\begin{Proposition}\label{Prop-FU-Damping2} There exists a numerical constant ${0<\bf c}$ such that if  $2 Gr \leq Re$ then we have the inequality $\ds{F \leq {\bf c} \left( \frac{1}{Gr} + 1\right)\, \frac{U^2}{\ell_0}}$. 
\end{Proposition}	
\pv  We consider here the differential equation  in (\ref{Eq_damped_NS}) where we must multiply each term by $\vf$ and integrate in the spatial variable. We recall that since the external force $\vf \in L^2(\Rt)$ satisfies the frequency localization given in (\ref{Loc_frec_force}) then it belongs to all the Sobolev spaces $H^s(\Rt)$ with $s \geq 0$ and we can write  
\begin{equation*}
\begin{split}
\int_{\Rt} \partial_t \vu(t,x) \cdot \vf(x) \, dx =&\, \nu \int_{\Rt} \Delta  \vu(t,x) \cdot \vf(x) \, dx - \int_{\Rt} \P((\vu\cdot \vec{\nabla}) \vu)(t,x)  \cdot \vf(x) \, dx + \Vert \vf \Vert^{2}_{L^2}\\
&\, - \beta  \int_{\Rt} \vu(t,x) \cdot \vf(x)\, dx. 
\end{split}
\end{equation*}
Then, we have 
\begin{equation}\label{Eq02}
\begin{split}
\Vert \vf \Vert^{2}_{L^2}  =& \int_{\Rt} \partial_t \vu(t,x) \cdot \vf(x) \, dx - \nu \int_{\Rt} \Delta \vu(t,x) \cdot \vf(x) \, dx + \int_{\Rt} \P((\vu\cdot \vec{\nabla}) \vu)(t,x)  \cdot \vf(x) \, dx \\
&\, + \beta \int_{\Rt} \vu(t,x) \cdot \vf(x)\, dx,
\end{split}
\end{equation}
where we must study each term on the right-hand side.  For the first term, as $\vf$ is a time-independent function we have
\begin{equation*}
\int_{\Rt} \partial_t \vu(t,x) \cdot \vf(x) \, dx = \partial_{t} \int_{\Rt} \vu(t,x)\cdot \vf(x) dx. 
\end{equation*}
For the second term, by the Cauchy-Schwarz inequality we obtain
\begin{equation*}
-\nu \int_{\Rt} \Delta \vu(t,x) \cdot \vf(x) \, dx = - \nu \int_{\Rt} \vu(t,x)\cdot \Delta  \vf(x) dx \leq \nu \Vert \vu(t,\cdot)\Vert_{L^2}\, \Vert \Delta  \vf \Vert_{L^2}.
\end{equation*}
For the third term, as $\vf$ is divergence-free vector field, by integrating by parts and by the H\"older inequalities, we write
\begin{equation*}
\begin{split}
\int_{\Rt} \P((\vu\cdot \vec{\nabla}) \vu)(t,x)  \cdot \vf(x) \, dx =&\, \int_{\Rt} (\vu \cdot \vec{\nabla}) \vu (t,x) \cdot \vf(x) dx = - \sum_{i,j=1}^{3} \int_{\Rt} u_i\, u_j (t,x) \partial_{j} f_i (x) dx \\
\leq &\, \Vert \vu(t,\cdot)\Vert^{2}_{L^2}\, \Vert \vec{\nabla} \otimes \vf \Vert_{L^\infty}. 
\end{split}
\end{equation*}
Finally, for the fourth term, by the Cauchy-Schwarz inequality we have
\begin{equation*}
\beta \int_{\Rt} \vu(t,x)\cdot \vf(x) dx \leq \beta \Vert \vu(t,\cdot)\Vert_{L^2}\, \Vert \vf \Vert_{L^2}. 
\end{equation*}
With these estimates at hand, we get back to the identity (\ref{Eq02}) to write
\begin{equation*}
\begin{split}
\Vert \vf \Vert^{2}_{L^2} \leq &\,\, \partial_{t} \int_{\Rt} \vu(t,x)\cdot \vf(x) dx + \nu \Vert \vu(t,\cdot)\Vert_{L^2}\, \Vert \Delta \vf \Vert_{L^2} + \Vert \vu(t,\cdot)\Vert^{2}_{L^2}\, \Vert \vec{\nabla} \otimes \vf \Vert_{L^\infty}  \\
&\, +  \beta  \Vert \vu(t,\cdot)\Vert_{L^2}\, \Vert \vf \Vert_{L^2}.
\end{split}
\end{equation*}
For $T>0$, in each term of this inequality we take the time average $\frac{1}{T} \int_{0}^{T}(\cdot) dt$; and always by the fact that $\vf$ does not depend on the time variable, by the Cauchy-Schwarz inequality (in the time variable) we get
\begin{equation*}
\begin{split}
\Vert \vf \Vert^{2}_{L^2} \leq &\,\, \frac{1}{T} \left( \int_{\Rt} \vu(T,x)\cdot \vf(x) dx - \int_{\Rt} \vu_0(x)\cdot \vf(x) dx \right) + \nu  \left( \frac{1}{T} \int_{0}^{T}  \Vert \vu(t,\cdot)\Vert^{2}_{L^2} dt \right)^{\frac12} \,  \Vert \Delta  \vf \Vert_{L^2} \\
&\,  + \left( \frac{1}{T}\,  \int_{0}^{T} \Vert \vu(t,\cdot)\Vert^{2}_{L^2} dt \right) \, \Vert \vec{\nabla} \otimes \vf \Vert_{L^\infty}   +  \beta  \left( \frac{1}{T} \int_{0}^{T} \Vert \vu(t,\cdot)\Vert^{2}_{L^2} dt \right)^{\frac12} \, \Vert \vf \Vert_{L^2}.
\end{split}
\end{equation*}
To study the first term on the right-hand side, we shall need the following estimate which is given by the damping parameter $-\beta \vu$ in the equations (\ref{Eq_damped_NS}).
\begin{Lemme} With the same hypothesis of Theorem \ref{Th1-Existence-NS-damped}, the solutions $\vu \in (L^{\infty}_{t})_{loc}L^{2}_{x} \cap (L^{2}_{t})_{loc}\dot{H}^{1}_{x}$ of the equation  (\ref{Eq_damped_NS}) verify the following control in time for all $t \leq 0$:
\begin{equation*}
\Vert \vu(t,\cdot)\Vert^{2}_{L^2} \leq e^{-\beta t }\Vert \vu_0 \Vert^{2}_{L^2}+ \frac{\Vert \vf \Vert^{2}_{\dot{H}^{-1}}}{\nu \beta}. 
\end{equation*}	
\end{Lemme}	
{\bf Proof.}  Consider the functions $\vu_{\delta_n}$ which are solutions of the regularized equation (\ref{NS-damped-regularized}). Our starting point is then the equality (\ref{energ_eq_delta}):
$$ \frac{d}{dt} \Vert \vu_{\delta_n}(t,\cdot)\Vert^{2}_{L^2} =  -2\nu \Vert \vu_{\delta_n}(t,\cdot)\Vert^{2}_{\dot{H}^1}+2\langle \vf,\vu_{\delta_n}(t,\cdot)\rangle_{\dot{H}^{-1}\times \dot{H}^1}-2\beta \Vert \vu_{\delta_n}(t,\cdot)\Vert^{2}_{L^2}.$$ 
Since  $\vf\in \dot{H}^{-1}(\Rt)$ we have $ \ds{2\langle \vf,\vu_{\delta_n}(t,\cdot)\rangle_{\dot{H}^{-1}\times \dot{H}^1}\leq \frac{\Vert \vf \Vert^{2}_{\dot{H}^{-1}}}{\nu}+\nu \Vert \vu_{\delta_n}(t,\cdot)\Vert^{2}_{\dot{H}^1}}$  and then we get 
\begin{equation*}
\frac{d}{dt} \Vert \vu_{\delta_n}(t,\cdot)\Vert^{2}_{L^2}\leq\frac{\Vert \vf\Vert^{2}_{\dot{H}^{-1}}}{\nu}  -\nu\Vert \vu_{\delta_n}(t,\cdot)\Vert^{2}_{\dot{H}^1} -2\beta \Vert \vu_{\delta_n}(t,\cdot)\Vert^{2}_{L^2}.
\end{equation*}
Then, we write 
$$ \frac{d}{dt} \Vert \vu_{\delta_n}(t,\cdot)\Vert^{2}_{L^2}\leq\frac{\Vert \vf\Vert^{2}_{\dot{H}^{-1}}}{\nu} - \beta \Vert \vu_{\delta_n}(t,\cdot)\Vert^{2}_{L^2},$$  and by the Gr\"onwall inequality we have the control:
$$\Vert \vu_{\delta_n}(t,\cdot)\Vert^{2}_{L^2}\leq e^{-\beta t}\Vert \vu_0 \Vert^{2}_{L^2}+\frac{\Vert \vf \Vert^{2}_{\dot{H}^{-1}}}{\nu \beta} \left( 1-e^{-\beta t}\right) \leq  e^{-\beta t}\Vert \vu_0 \Vert^{2}_{L^2}+\frac{\Vert \vf \Vert^{2}_{\dot{H}^{-1}}}{\nu \beta},$$ 
for all time $t\in [0,+\infty[$. Now, we will  recover this control in time  for the limit function $\vu$: we regularize in the time variable the quantity $\Vert \vu_{\delta_n}(t,\cdot)\Vert^{2}_{L^2}$ by  a convolution product with a positive   function $w\in\mathcal{C}^{\infty}_{0}([-\eta,\eta])$ (for $\eta>0$) such that $\ds{\int_{\mathbb{R}}w(t)dt=1}$.  In this way, in the previous inequality we have $$\Vert w\ast \vu_{\delta_n}(t,\cdot)\Vert_{L^2}\leq  w\ast \Vert \vu_{\delta_n}(t,\cdot)\Vert^{2}_{L^2}\leq w\ast \left( e^{-\beta t}\Vert \vu_0 \Vert^{2}_{L^2}+\frac{\Vert \vf \Vert^{2}_{\dot{H}^{-1}}}{\nu \beta}\right).$$ 
Moreover, since $(\vu_{\delta_n})_{n\in\mathbb{N}}$ converges  weakly$-*$ to $\vu$ in $(L^{\infty}_{t})_{loc}(L^{2}_{x})$ then  $w\ast \vu_{\delta_n}(t,\cdot)$ converges weakly$-*$  to $w\ast\vu(t,\cdot)$ in  $L^2(\Rt)$ and then we can write    
\begin{equation*}
\Vert w\ast \vu(t,\cdot)\Vert^{2}_{L^2} \leq  \liminf_{n\longrightarrow +\infty} \Vert w\ast \vu_{\delta_n}(t,\cdot)\Vert^{2}_{L^2}\leq  w\ast \left( e^{-\beta t}\Vert \vu_0 \Vert^{2}_{L^2}+\frac{\Vert \vf \Vert^{2}_{\dot{H}^{-1}}}{\nu \beta}\right).
\end{equation*} Finally, for $t$ a Lebesgue point of the function $t\mapsto \Vert \vu(t,\cdot)\Vert^{2}_{L^2}$  we have  the wished  control in time and we can extend this inequality to all time $t\in [0,+\infty[$ by the weak continuity of the function $t\mapsto \Vert \vu(t,\cdot)\Vert^{2}_{L^2}$.  \finpv

\medskip

Once we have this estimate at our disposal, in the first term on right-hand side,  first we apply the Cauchy-Schwarz inequality and by this estimate we obtain 
\begin{equation*}
\begin{split}
\Vert \vf \Vert^{2}_{L^2} \leq &\,\, \frac{1}{T} \left(e^{-\beta T} \Vert \vu_0\Vert^{2}_{L^2} + \frac{1}{\nu \beta}\Vert \vf \Vert^{2}_{\dot{H}^{-1}}+ \Vert \vu_0 \Vert^{2}_{L^2} + \Vert \vf \Vert^{2}_{L^2}  \right) + \nu  \left( \frac{1}{T} \int_{0}^{T} \Vert \vu(t,\cdot)\Vert^{2}_{L^2} dt \right)^{\frac12} \,  \Vert  \Delta  \vf \Vert_{L^2} \\
&\,  + \left( \frac{1}{T}\,  \int_{0}^{T} \Vert \vu(t,\cdot)\Vert^{2}_{L^2} dt \right) \, \Vert \vec{\nabla} \otimes \vf \Vert_{L^\infty}   +  \beta  \left( \frac{1}{T} \int_{0}^{T} \Vert \vu(t,\cdot)\Vert^{2}_{L^2} dt \right)^{\frac12} \, \Vert \vf \Vert_{L^2}.
\end{split}
\end{equation*}
In each term of this inequality, we divide by $\ell^{3}_{0}$ and by $T$, we also apply the $\limsup$ when $T \to +\infty$. Thus, by definition of the quantities $F$ and $U$ given in (\ref{Def_F}) and (\ref{Def_U_Intro}) respectively, we get
\begin{equation*}
F^2 \leq \nu \,  U \frac{\Vert \Delta  \vf \Vert_{L^2} }{\ell^{3/2}_{0}}+ U^2\, \Vert \vec{\nabla} \otimes \vf \Vert_{L^\infty} + \beta \, UF, 
\end{equation*}
hence, 
\begin{equation}\label{Eq06}
F \leq \nu \,  U \frac{\Vert \Delta \vf \Vert_{L^2} }{\ell^{3/2}_{0}\, F}+ U^2 \frac{\Vert \vec{\nabla} \otimes \vf \Vert_{L^\infty}}{F} + \beta \, U, 
\end{equation}
where we still need to study the first term and the second term on the right-hand side.  For the first term,  we recall that $\ell^{3/2}_{0} F = \Vert \vf \Vert_{L^2}$. Moreover, by the  frequency localization of $\vf$ given in (\ref{Loc_frec_force}),  for the first term we have
\begin{equation}\label{Eq07}
\frac{\Vert \Delta \vf \Vert_{L^2} }{\ell^{3/2}_{0}\, F} =  \frac{\Vert \Delta  \vf \Vert_{L^2} }{\Vert \vf \vert_{L^2	}} \leq \mathfrak{c}^2\, \frac{\Vert \vf \Vert_{L^2}}{\ell^{2}_{0} \Vert \vf \Vert_{L^2}} = \frac{\mathfrak{c}^2}{\ell^{2}_{0}}, 
\end{equation} 
where $0<\mathfrak{c}$ is the dimensionless constant in (\ref{Loc_frec_force}). 
For the second term,  always by  (\ref{Loc_frec_force}) and by the  Bernstein  inequalities (see \cite[Chapter $2$]{Grafakos}) there exists a numerical constant $0<c_0$  such that we have
\begin{equation*}
\frac{\Vert \vec{\nabla} \otimes \vf \Vert_{L^\infty}}{F} \leq \frac{c_0}{\ell_0} 
\end{equation*}
Thus, we set a constant ${0<\bf c}$ such that $\max\left( \mathfrak{c}^{2}, c_0,1\right) \leq {\bf c}$, and by these estimates we can write
\begin{equation*}
F \leq {\bf c} \left( \frac{\nu \, U}{\ell^{2}_{0}} + \frac{U^2}{\ell_0}\right)+ \beta \, U.
\end{equation*} 
Thereafter,  since   $\ds{\beta= \frac{\nu}{\ell^{2}_{0}}}$,  by recalling that $Re= \frac{U \ell_0}{\nu}$, and moreover, if we assume that $2Gr \leq Re$ then   we can  write 
\begin{equation*}
F \leq {\bf c} \left( \frac{2\nu \, U}{\ell^{2}_{0}} + \frac{U^2}{\ell_0}\right) =  {\bf c}  \left( \frac{2\nu}{U\, \ell_0} + 1 \right) \frac{U^2}{\ell_0}   =  {\bf c}  \left( \frac{2}{Re} + 1 \right) \frac{U^2}{\ell_0}\leq {\bf c}  \left( \frac{1}{Gr} + 1 \right) \frac{U^2}{\ell_0}. 
\end{equation*} 
\finpv
\begin{Proposition}\label{Prop-Epsilon-F-U-Damping2} For the quantities $\mathcal{E}$, $U$ and $F$  defined in (\ref{Def_Epsilon_Intro}), (\ref{Def_U_Intro}) and (\ref{Def_F}) respectively, we have the following inequality  $\mathcal{E} \leq FU$.
\end{Proposition}	
\pv By the energy inequality (\ref{Energy-inequality}) and by the Cauchy-Schwarz inequality (first in the spatial variable and then in the temporal variable) we can write 
\begin{equation*}
\begin{split}
2 \nu \int_{0}^{T} \Vert  \vu(t,\cdot) \Vert^{2}_{\dot{H}^1}\, dt \leq &\,  \Vert \vu_0 \Vert^{2}_{L^2} + 2 \int_{0}^{T} \int_{\Rt} \vf(x)\cdot \vu(t,x) \, dx\, dt  \\
\leq &\, \Vert \vu_0 \Vert^{2}_{L^2} + 2 \int_{0}^{T} \Vert \vf \Vert_{L^2}\, \Vert \vu(t,\cdot)\Vert_{L^2}\, dt  \\
\leq & \, \Vert \vu_0 \Vert^{2}_{L^2}  +  \Vert \vf \Vert_{L^2}\,  T^{1/2}\, \left( \int_{0}^{T} \Vert \vu(t,\cdot) \Vert^{2}_{L^2}\, dt\right)^{\frac12}.
\end{split}
\end{equation*}
Thereafter, in each term of this inequality we divide by $\ell^{3}_{0}$ and by $T$; and we also apply the $\limsup$ when $T \to +\infty$. By definition of the quantities $\mathcal{E}$, $F$ and $U$ we obtain the wished inequality.  \finpv


Now we are able to prove the estimate (\ref{Estimate_above_Kolmogorov_Damping2}). In the inequality  given in Proposition \ref{Prop-FU-Damping2} we multiply by $U$ to get $\ds{FU \leq {\bf c} \left(\frac{1}{Gr} + 1 \right) \frac{U^3}{\ell_0}}$. Then, by the inequality given in Proposition \ref{Prop-Epsilon-F-U-Damping2} we obtain the wished estimate (\ref{Estimate_above_Kolmogorov_Damping2}).  Theorem \ref{Th_Kolmogorow_Law_Damping2} is proven. \finpv


\end{document}